\documentclass[11pt, reqno]{amsart}    	
\usepackage{geometry}                		
\usepackage{graphicx}				
\usepackage{amssymb}
\usepackage{amsmath}
\usepackage{comment}
\usepackage{amsthm}
\usepackage{mathrsfs}
\allowdisplaybreaks
\newtheorem{theorem}{Theorem}
\numberwithin{theorem}{section}
\newtheorem{proposition}[theorem]{Proposition}
\newtheorem{lemma}[theorem]{Lemma}
\newtheorem{corollary}[theorem]{Corollary}
\newcommand{\EndProof}{\hfill $\square$}
\usepackage[colorlinks=true, linkcolor=blue,
 citecolor=blue, urlcolor=blue]{hyperref}
 \usepackage{dsfont}
 \usepackage{tikz-cd}
 \usepackage{caption}
 \usepackage{enumitem} \usepackage{mathdots}

\title[The Twisted Satake Transform and the Casselman-Shalika Formula]{The Twisted Satake Transform and the Casselman-Shalika Formula for Quasi-Split Groups}
\author{Nadya Gurevich and Edmund Karasiewicz}

\address{Nadya Gurevich: Department of Mathematics, Ben Gurion University of the Negev, Be'er Sheva,  Israel 8410501}
\email{ngur@math.bgu.ac.il }
\address{Edmund Karasiewicz: Department of Mathematics, Ben Gurion University of the Negev, Be'er Sheva,  Israel 8410501}
\email{karasiew@post.bgu.ac.il}

\begin{document}
\maketitle 
\begin{abstract}
We prove the Casselman-Shalika formula for unramified groups over a non-archimedean local field by studying the action of the spherical Hecke algebra on the space of compact spherical Whittaker functions via the twisted Satake transform. This method provides a conceptual explanation of the appearance of characters of a dual group in the Casselman-Shalika formula.
\end{abstract}

\section{Introduction}

In this article we provide a conceptual explanation of the appearance of characters of a complex reductive group in the Casselman-Shalika formula of a connected unramified group $G$ over a nonarchimedean local field $F$. This extends the work of the first author in \cite{N13}, where $G$ was assumed to be split of adjoint type. 

For this introduction we assume that $G$ is split unless otherwise specified. Let  $T$ be a maximal torus, $K$ a maximal compact subgroup, and $U$ a maximal unipotent subgroup of $G$ such that $G=UTK$ (Iwasawa Decomposition). Let $\Psi:U\rightarrow \mathbb{C}^{\times}$ be a non-degenerate character of conductor $\mathfrak{p}$. (For definitions see Subsection \ref{WF} and Section \ref{TwistedSatakeSec}. Casselman-Shalika \cite{CS80} assume $\Psi$ is trivial on $U\cap K$ and nontrivial on any subgroup with a larger abelianization. They call such a $\Psi$ unramified; we will say it has conductor $\mathcal{O}$.) Let $\mathcal{H}_{K}\stackrel{\mathrm{def}}{=}C_{c}^{\infty}(K\backslash G/K)$ be the spherical Hecke algebra. The algebra $\mathcal{H}_{K}$ acts on $\mathrm{Ind}_{U}^{G}(\Psi)^{K}=\{f:G\rightarrow \mathbb{C}|f(ugk)=\Psi(u)f(g),\text{ for all }g\in G,u\in U,k\in K\}$ by right convolution. For $h\in \mathcal{H}_{K}$ we write $\hat{h}(g)=h(g^{-1})$. A spherical Whittaker function on $G$ with respect to $K$, $\Psi$, and a character $\chi:\mathcal{H}_{K}\rightarrow \mathbb{C}$ is a function $\mathcal{W}_{\chi}:G\rightarrow \mathbb{C}$ that satisfies the following properties.
\begin{enumerate}
\item $\mathcal{W}_{\chi}\in\mathrm{Ind}_{U}^{G}(\Psi)^{K}$.\label{WhitTrans}
\item $\mathcal{W}_{\chi}*h=\chi(\hat{h})\mathcal{W}_{\chi}$, for all $h\in \mathcal{H}_{K}$.\label{Eigen}
\end{enumerate}
Since $G=UTK$, item (\ref{WhitTrans}) implies that $\mathcal{W}_{\chi}$ is determined by its values on $T/T\cap K$; item (\ref{Eigen}) in conjunction with the multiplicity one theorem for Whittaker models implies that the associated eigenvalues determine $\mathcal{W}_{\chi}$ up to scaling. The Casselman-Shalika formula (\cite{CS80}, Theorem 5.4) is a formula for the function $\mathcal{W}_{\chi}$ evaluated at the elements of $T/T\cap K$. 

Casselman-Shalika's original approach \cite{CS80} is based on decomposing the spherical Whittaker function using Casselman's basis. One function in this decomposition can be computed directly; the others can then be computed via a symmetry. One can manipulate the resulting formula and apply the Weyl character formula to find a character of a finite dimensional representation of the Langlands dual group of $G$. When $G$ is not split one can make a similar statement. In this case the description of the dual group (not necessarily the Langlands dual group) seems to first appear in Tamir \cite{T91}. 

While it is not difficult to identify this character as a byproduct of the method of Casselman-Shalika, Gurevich \cite{N13} provided a conceptual explanation of its appearance when $G$ is split and adjoint. Gurevich's proof is based on studying the $\mathcal{H}_{K}$-module $\mathrm{ind}_{U}^{G}(\overline{\Psi})^{K}=\{f\in\mathrm{Ind}_{U}^{G}(\overline{\Psi})^{K}|\mathrm{supp}(f)\text{ is compact mod }U\}$, where $\overline{\Psi}$ is the complex conjugate of $\Psi$. The character naturally appears through the Satake isomorphism.

In addition to providing a conceptual explanation of the appearance of the character, Gurevich's proof also avoids the multiplicity one theorem for Whittaker functionals, which is invoked by Casselman-Shalika \cite{CS80}. Thus this method may provide an alternative approach to the Casselman-Shalika formula in the the case of covering groups, first proved by Chinta-Offen \cite{CO13} and McNamara \cite{M16}, and $p$-adic loop groups, first proved by Patnaik \cite{P17}. 

However, Gurevich's proof does not determine the normalization factor $\zeta(\chi)$ (Casselman-Shalika \cite{CS80}, Theorem 5.4), which characterizes the spherical representations that do not possess a Whittaker functional. Nevertheless, almost all spherical representations are principal series, which do possess a Whittaker functional. In these cases the choice of a normalization of the Whittaker functional determines the constant of proportionality in our formula.

The present paper extends Gurevich's proof of the Casselman-Shalika formula to connected unramified groups, the class of groups considered by Casselman-Shalika \cite{CS80}. 

The foundation for our proof is an explicit $\mathcal{H}_{K}$-module isomorphism. To state this precisely we introduce some additional notation. The torus $T$ has a cocharacter lattice $X_{*}(T)$ and Weyl group $W=N_{G}(T)/T$. Using the isomorphism $T/T\cap K\cong X_{*}(T)$, for any $\mu\in X_{*}(T)$ we fix an element $t_{\mu}\in T$ representing $\mu$ in the quotient.

We will write $\mathbb{C}[X_{*}(T)]$ for the group ring of $X_{*}(T)$ and let $\{e^{\mu}|\mu\in X_{*}(T)\}$ be the natural basis. The action of $W$ on $T$ induces an action on $\mathbb{C}[X_{*}(T)]$ and we write $\mathbb{C}[X_{*}(T)]^{W}$ and $\mathrm{alt}(\mathbb{C}[X_{*}(T)])$ for the set of symmetric and alternating elements, respectively. 

Let $B=TU\subset G$ be a Borel subgroup with modular character $\delta$  and unipotent radical $U$. Let $X_{*}(T)^{+}$ $(X_{*}(T)^{++})$ be the set of dominant (strictly dominant) cocharacters with respect to $B$. Let $\rho^{\vee}$ be one half the sum of the dominant coroots.

Let $\widehat{G}$ be the Langlands dual group of $G$ and $\mathrm{Rep}(\widehat{G})\cong \mathbb{C}[X_{*}(T)]^{W}$ its $\mathbb{C}$-algebra of finite dimensional characters. The set of characters of highest weight representations $\{\mathrm{ch}V_{\lambda}| \lambda\in X_{*}(T)^{+}\}$ is a basis for $\mathrm{Rep}(\widehat{G})$. Using the Satake isomorphism $\mathcal{S}:\mathcal{H}_{K}\rightarrow \mathbb{C}[X_{*}(T)]^{W}$, for every $\lambda\in X_{*}(T)^{+}$ we can define $A_{\lambda}=\mathcal{S}^{-1}(\mathrm{ch}V_{\lambda})\in \mathcal{H}_{K}$.

The space $\mathrm{ind}_{U}^{G}(\overline{\Psi})^{K}$ has a basis $\{\phi_{\mu}|\mu\in X_{*}(T)^{++}\}$ defined as follows. Consider the $G$-invariant pairing $\langle-,-\rangle_{\Psi}:\mathrm{ind}_{U}^{G}(\overline{\Psi})\otimes \mathrm{Ind}_{U}^{G}(\Psi)\rightarrow \mathbb{C}$ defined by $\langle h,f\rangle_{\Psi}=\int_{U\backslash G}h(g)f(g)dg$. For each $\mu\in X_{*}(T)^{++}$ the function $\phi_{\mu}\in \mathrm{ind}_{U}^{G}(\overline{\Psi})^{K}$ is defined by the property that 
\begin{equation}\label{evalIntro}
\delta(t_{\mu})^{-1/2}f(t_{\mu})=\langle\phi_{\mu},f\rangle_{\Psi},
\end{equation}
for all $f\in \mathrm{Ind}_{\overline{U}}^{G}(\Psi)^{K}$.

Now we can state the main technical theorem.
\begin{theorem}\label{HKIsoIntro}
The linear map $j:\mathrm{ind}_{U}^{G}(\Psi)^{K}\rightarrow \mathrm{alt}(\mathbb{C}[X_{*}(T)])$ defined by $j(\phi_{\mu})=\mathrm{alt}(e^{\mu})$ for $\mu\in X_{*}(\textbf{A})^{++}$ is an $\mathcal{H}_{K}\stackrel{\mathcal{S}}{\cong} \mathbb{C}[X_{*}(T)]^{W}$-module isomorphism.
\end{theorem}

When $G$ is adjoint then $\rho^{\vee}\in X_{*}(T)^{++}$, and Theorem \ref{HKIsoIntro} quickly yields the Casselman-Shalika formula. Specifically, the Weyl character formula and Theorem \ref{HKIsoIntro} imply for $\mu\in X_{*}(T)^{++}$ and $\lambda \in X_{*}(T)^{+}$
\begin{equation}\label{RhoIntro}
j(\phi_{\rho^{\vee}}*A_{\lambda})=j(\phi_{\rho^{\vee}})\cdot \mathrm{ch}V_{\lambda}=\mathrm{alt}(e^{\rho^{\vee}})\cdot \mathrm{ch}V_{\lambda}=\mathrm{alt}(e^{\lambda+\rho^{\vee}})=j(\phi_{\lambda+\rho^{\vee}}).
\end{equation}
Thus $\phi_{\rho^{\vee}}*A_{\lambda}=\phi_{\lambda+\rho^{\vee}}$. Therefore because $\mathcal{W}_{\chi}$ is an $\mathcal{H}_{K}$-eigenfunction and the pairing is $G$-invariant, for any $\lambda\in X_{*}(T)^{+}$ we have 
\begin{multline}\label{RecurIntro}
\chi(A_{\lambda})\delta(t_{\rho^{\vee}})^{-1/2}\mathcal{W}_{\chi}(t_{\rho^{\vee}})=\langle\phi_{\rho^{\vee}},\mathcal{W}_{\chi}*\hat{A}_{\lambda}\rangle_{\Psi}=\langle \phi_{\rho^{\vee}}*A_{\lambda},\mathcal{W}_{\chi}\rangle_{\Psi}\\=\langle \phi_{\lambda+\rho^{\vee}},\mathcal{W}_{\chi}\rangle_{\Psi}=\delta(t_{\lambda+\rho^{\vee}})^{-1/2}\mathcal{W}_{\chi}(t_{\lambda+\rho^{\vee}}).
\end{multline}

\begin{theorem}[Casselman-Shalika Formula for Adjoint Groups]\label{CSAdIntro} For any $\lambda\in X_{*}(T)^{+}$
\begin{equation}
\mathcal{W}_{\chi}(t_{\lambda+\rho^{\vee}})=\delta^{1/2}(t_{\lambda})\chi(\mathcal{S}^{-1}(\mathrm{ch}V_{\lambda}))\mathcal{W}_{\chi}(t_{\rho^{\vee}}).
\end{equation}
\end{theorem}
Since $G$ is adjoint Theorem \ref{CSAdIntro} directly implies a formula for characters $\Psi$ of conductor $\mathcal{O}$, which is important for global applications. For details see Section \ref{CSAdjointSec}.

If $G$ is not adjoint, then it may be that $\rho^{\vee}\notin X_{*}(T)$, in which case (\ref{RhoIntro}) is meaningless. Nevertheless we prove a substitute (Lemma \ref{TensorDecomp}). For any $\lambda\in X_{*}(T)^{+}$ and $\mu\in X_{*}(T)^{++}$ define $c_{\mu,\lambda}^{\eta}\in\mathbb{C}$ such that $\mathrm{ch}V_{\mu-\rho^{\vee}}\cdot\mathrm{ch}V_{\lambda}=\sum_{\eta}c_{\mu,\lambda}^{\eta}\mathrm{ch}V_{\eta-\rho^{\vee}}$. (In case $\rho^{\vee}\notin X_{*}(T)$ we view $\mathrm{ch}V_{\mu-\rho^{\vee}}$ as a representation of the simply connected cover of $\widehat{G}$.) Then

\begin{equation}\label{TensorDecompEqnIntro}
\phi_{\mu}*A_{\lambda}=\sum_{\eta\in X_{*}(\textbf{A})^{++}}c_{\mu,\lambda}^{\eta}\phi_{\eta}.
\end{equation}
Lines (\ref{evalIntro}) and (\ref{TensorDecompEqnIntro}) yield a family of recursions for $\mathcal{W}_{\chi}$. For any $\lambda\in X_{*}(T)^{+}$ and $\mu\in X_{*}(T)^{++}$
\begin{equation}\label{WhittakerRecurIntro}
\delta^{-1/2}(t_{\mu})\chi(\mathcal{S}^{-1}(\mathrm{ch}V_{\lambda}))\,\mathcal{W}_{\chi}(t_{\mu})=\sum_{\eta\in X_{*}(T)^{++}}c_{\mu,\lambda}^{\eta}\delta^{-1/2}(m_{\eta})\mathcal{W}_{\chi}(t_{\eta}).
\end{equation}

In Proposition \ref{dimV} we show that these recursions have a solution space of dimension at most one. Moreover, a nonzero solution to these recursions is provided by characters of the simply connected cover of $\widehat{G}$, for almost all $\chi$. This yields the Casselman-Shalika formula (Theorem \ref{CSFormulaFinal}).

So far we have focused on the case where $G$ is split. This was a matter of convenience and in the main body of this paper we will work with an arbitrary connected unramified group $G$. The most important new feature required for this more general case is identifying the correct dual group ${^{L}G^{\dagger}}$ (Subsection \ref{DualGroup}). In fact, we work in the context of the universal principal series, which highlights the role played by this dual group.

The above proof of the Casselman-Shalika formula only works for characters of conductor $\mathfrak{p}$. However, for global applications one must consider characters of conductor $\mathcal{O}$. We include an alternative proof of the Casselman-Shalika formula for characters of conductor $\mathfrak{p}$ via a reduction to the adjoint case. The advantage of this approach is that it allows us to recover the Casselman-Shalika formula for characters of conductor $\mathcal{O}$.

We conclude this introduction with an outline of this paper. In Section \ref{Notation} we establish notation and recall relevant theorems. In Section \ref{IwahoriHecke} we recall the Bernstein presentation for the Iwahori-Hecke algebra of a connected unramified group (Theorem \ref{BP}). This presentation is used in Section \ref{SavinsIsoSec} where we generalize Savin's isomorphism \cite{S12} to connected unramified groups (Theorem \ref{SavinIso}).

In Section \ref{TwistedSatakeSec} we determine the $\mathcal{H}_{K}$-module structure of $\mathrm{ind}_{\overline{U}}^{G}(\overline{\Psi})^{K}$ through a study of the twisted Satake transform $\mathcal{S}_{\Psi}$. The main technical result is the computation of the kernel of $\mathcal{S}_{\Psi}$ (Lemma \ref{TwistedKer}). It is in this computation where we use that $\Psi$ has conductor $\mathfrak{p}$.

In Section \ref{UWF} we study a spherical Whittaker function of conductor $\mathfrak{p}$. The main result of this section is a recursive formula for this function (Proposition \ref{WhittakerRecur}). In Section \ref{CSAdjointSec} we specialize this to the case where $G$ is adjoint and prove the Casselman-Shalika formula for characters of conductor $\mathfrak{p}$ (Theorem \ref{CSAdjoint}) and $\mathcal{O}$ (Proposition \ref{CSAdO}). In Section \ref{CSFormulaSec} we we generalize the argument from Section \ref{CSAdjointSec} to connected unramified groups. Unfortunately, this does not yield any information about the conductor $\mathcal{O}$ case. Therefore, in Section \ref{CSGeneral} we prove the Casselman-Shalika formula for connected unramified groups and characters of conductor $\mathfrak{p}$ via a reduction to the adjoint case. This approach allows us to recover the Casselman-Shalika formula for characters of conductor $\mathcal{O}$, which is described in Section \ref{ConductorO}.

\section{Notation}\label{Notation}
\subsection{Fields}

Let $F$ be a nonarchimedean local field with finite residue field. If $L$ is a field extension of $F$, we write $\mathcal{O}_{L}$ for the ring of integers of $L$ with maximal ideal $\mathfrak{p}_{L}$, and $\kappa_{L}$ for the residue field. If $\kappa_{L}$ is finite let $q_{L}=|\kappa_{L}|$. When $L=F$ we may suppress the subscripts. Let $\mathrm{ord}$ be the discrete valuation of $F$ with value group $\mathbb{Z}$. Let $E$ be a finite unramified extension of $F$. Since $E\supseteq F$ is unramified the canonical extension of $\mathrm{ord}$ from $F$ to $E$ also has value group $\mathbb{Z}$. 

\subsection{Algebraic Groups}\label{AlgGrps}

Throughout this paper we use boldface characters for group schemes over $F$, such as \textbf{H}, and plain text characters for their group of $F$-points, such as $H$. 

Let \textbf{H} be a group scheme defined over $F$. We write $Z(\textbf{H})$ for the center of $\textbf{H}$ and $\mathcal{D}\textbf{H}$ for the derived subgroup of $\textbf{H}$. For any field $L\supseteq F$, let $X^{*}(\textbf{H})_{L}$ be the group of algebraic characters defined over $L$ and $X_{*}(\textbf{H})_{L}$  the group of algebraic cocharacters defined over $L$. If $L=F_{s}$, a separable closure of $F$, then we omit the subscript.

Let $\textbf{G}$ be a connected quasi-split reductive group scheme defined over $F$ that is split over $E$. (i.e. $\textbf{G}$ is unramified.) Inside of $\textbf{G}$ we fix a minimal parabolic subgroup $\textbf{P}$ containing a maximal $F$-split torus $\textbf{A}$ with centralizer $\textbf{M}$ and normalizer $\textbf{N}$. Since $\textbf{G}$ is quasi-split $\textbf{M}$ is a maximal torus in $\textbf{G}$. Thus $X_{*}(\textbf{M})_{F}\cong X_{*}(\textbf{A})$. The parabolic subgroup $\textbf{P}$ admits a Levi decomposition $\textbf{P}=\textbf{M}\textbf{U}$, where $\textbf{U}$ is the unipotent radical of $\textbf{P}$. We write $\overline{\textbf{P}}$ and $\overline{\textbf{U}}$ for the opposite parabolic and its unipotent radical respectively. 

The maximal split torus $\textbf{A}$ in $\textbf{G}$ determines a relative root datum $(X^{*}(\textbf{A}),\Phi,X_{*}(\textbf{A}),\Phi^{\vee})$, where $\Phi=\Phi(\textbf{G},\textbf{A})$ is the set of relative roots with respect to $\textbf{A}$. For each $\alpha\in \Phi$ we write $\textbf{U}_{\alpha}$ for the root subgroup associated to $\alpha$. (Note that if $2\alpha\in \Phi$, then $\textbf{U}_{2\alpha}\subset \textbf{U}_{\alpha}.$) Our choice of $\textbf{P}$ identifies a set of positive roots $\Phi^{+}$, from which we can extract a set of simple roots $\Delta$. The set of negative roots is $\Phi^{-}=-\Phi^{+}$. Let $\langle-,-\rangle: X^{*}(\textbf{A})\times X_{*}(\textbf{A})\rightarrow \mathbb{Z}$ be the pairing defined by  $\lambda\circ\mu(t)=t^{\langle\lambda,\mu\rangle}$. Let $\Lambda^{\vee}=\{v\in \mathrm{span}_{\mathbb{Z}}(\Phi^{\vee})\otimes \mathbb{Q}|\langle\alpha,v\rangle\in\mathbb{Z}\text{ for all }\alpha\in \Phi\}$ be the coweight lattice in $\mathrm{span}_{\mathbb{Z}}(\Phi^{\vee})\otimes \mathbb{Q}$ with respect to $\Phi$.

Let $W$ be the relative Weyl group associated to the root system $\Phi$, which is isomorphic to $N/M$. The action of $W$ on $\textbf{A}$ induces an action on $X_{*}(\textbf{A})$ which extends to an action on $\mathbb{C}[X_{*}(\textbf{A})]$. Define $\mathrm{alt}:\mathbb{C}[X_{*}(\textbf{A})]\rightarrow\mathbb{C}[X_{*}(\textbf{A})]$ by $\mathrm{alt}(f)=\sum_{w\in W}\mathrm{sign}(w)(w\cdot f)$, where $\mathrm{sign}$ is the sign character of $W$. Let $\mathrm{alt}(\mathbb{C}[X_{*}(\textbf{A})])$ denote the image of $\mathrm{alt}$.

The subgroup $M$ acts on $U$ by conjugation. The modular character of $P=MU$ is defined by $d(pup^{-1})=\delta_{P}(p)du$, where $du$ is a Haar measure of $U$. Similarly we defined $\delta_{\overline{P}}$ to be the modular character of $\overline{P}=M\overline{U}$. Note that $\delta_{P}=\delta_{\overline{P}}^{-1}$.

An element $\lambda\in X_{*}(\textbf{A})$ is dominant (strictly dominant) with respect to $P$ if $\langle\lambda,\alpha\rangle\geq 0$ ($\langle\lambda,\alpha\rangle> 0$) for all $\alpha\in \Delta$. Let $X_{*}(\textbf{A})^{+}$ ($X_{*}(\textbf{A})^{++}$) denote the set of dominant (strictly dominant) elements of $X_{*}(\textbf{A})$ with respect to $P$.

Later we will have three connected reductive groups $\textbf{G}$, $\textbf{G}^{\prime}$, and $\textbf{G}^{\prime\prime}$ appearing simultaneously. All the notation that we introduced for $\textbf{G}$ will be carried over to $\textbf{G}^{\prime}$ and $\textbf{G}^{\prime\prime}$ and augmented with one or two primes respectively.

\subsection{Bruhat-Tits Theory}

It will be convenient to establish some terminology from Bruhat-Tits theory. We follow Tits \cite{T79}. (For additional details see Bruhat-Tits \cite{BT72,BT84}, and Vigneras \cite{V16}, Section 3.) Let $\mathscr{B}$ denote the building of $G$ over $F$ and let $\mathscr{A}\subset \mathscr{B}$ be the (enlarged) apartment associated to the maximal $F$-split torus $A$. We fix a hyperspecial vertex $x_{0}\in\mathscr{A}$, which we use to make an identification $\mathscr{A}\cong X_{*}(\textbf{A})\otimes_{\mathbb{Z}} \mathbb{R}$.

The Bruhat-Tits homomorphism $\nu:M\rightarrow X_{*}(\textbf{M})_{F}\cong X_{*}(\textbf{A})$ is characterized by the equations $\langle\chi,\nu(m)\rangle=-\mathrm{ord}(\chi(m))$, for all $\chi\in X^{*}(\textbf{M})_{F}$. Let $M^{\circ}=\ker\nu$, then we have the exact sequence
\begin{equation}\label{ExactSq1}
1\rightarrow M^{\circ}\rightarrow M\stackrel{\nu}{\rightarrow} X_{*}(\textbf{A})\rightarrow 1.
\end{equation}
(For surjectivity see Cartier \cite{C79}, page 135.)

For convenience, we fix a splitting of this sequence $s:X_{*}(\textbf{A})\rightarrow M$ (i.e. $s$ is a group homomorphism such that $\nu\circ s=id_{X_{*}(\textbf{A})}$), which exists because $M$ is abelian and $X_{*}(\textbf{A})$ is a free abelian group, and write $m_{\lambda}=s(\lambda)$. 

Let $\alpha+k$ be an affine function on $X_{*}(\textbf{A})\otimes \mathbb{R}$, where $\alpha\in \Phi$ and $k\in \mathbb{R}$. The group $U_{\alpha}$ contains the subgroups $X_{\alpha+k}$ (Tits \cite{T79}, Section 1.4). These groups satisfy $X_{\alpha+k}\subseteq X_{\alpha+\ell}$ if and only if $k\geq \ell$ and determine a filtration of $U_{\alpha}$. This filtration can be used to define the set of affine roots $\Phi_{\mathrm{aff}}$ as in Tits \cite{T79}, Section 1.6.

Under the identification $\mathscr{A}\cong X_{*}(\textbf{A})\otimes \mathbb{R}$, the cone $X_{*}(\textbf{A})^{+}\otimes_{\mathbb{Z}_{\geq0}} \mathbb{R}_{\geq0}$ identifies a unique conical chamber $\mathscr{C}$ in $\mathscr{A}$ with apex $x_{0}$. Let $C$  be the unique chamber in $\mathscr{C}$ containing $x_{0}$ and define $I$ to be the Iwahori subgroup associated to $C$ (Tits \cite{T79}, Section 3.7). Recall that $I$ is contained in the fixator of $C$. The stabilizer of $x_{0}$ is a maximal compact subgroup $K=\mathrm{Stab}_{G}(x_{0})$. 

The next few propositions collect a few facts that will be useful in later computations.

\begin{proposition}\label{ConjContract} Let $\alpha\in \Phi$, $\lambda\in X_{*}(\textbf{A})$, and $k\in \mathbb{R}$. 
\begin{enumerate}
\item\label{conj}  $m_{\lambda}X_{\alpha+k}m_{\lambda}^{-1}=X_{\alpha+k-\alpha(\lambda)}$. 
\item\label{contract} Thus if $\lambda$ is dominant with respect to $P$ and $\alpha\in \Phi^{-}$, then $m_{\lambda}X_{\alpha+k}m_{\lambda}^{-1}\subseteq X_{\alpha+k}$.
\end{enumerate}
\end{proposition}
\textbf{Proof:} Item (\ref{conj}) is stated in Tits \cite{T79}, Section 1.4.2. Item (\ref{contract}) follows from item (\ref{conj}).\EndProof\\

If $\lambda$ is dominant we will informally say that $m_{\lambda}$ contracts $\overline{U}$. Furthermore, if $\lambda$ is strictly dominant, then this containment is strict and we will say $m_{\lambda}$ strictly contracts $\overline{U}$.

\begin{proposition}\label{FactsProp1}$ $
\begin{enumerate}
\item\label{Fact1} $M^{\circ}=M\cap K$
\item\label{Fact2} $P\cap K=(M\cap K)(U\cap K)$
\item\label{Fact3}\label{IndexModChar} Let $\lambda\in X_{*}(\textbf{A})^{+}$. Then $[Im_{\lambda}I:I]=\delta_{\overline{P}}(m_{\lambda})^{-1}$. 
\item\label{Iwasawa}$G=\overline{U}MK$. (Iwasawa Decomposition) 
\end{enumerate}
\end{proposition}

\textbf{Proof:} For item (\ref{Fact1}) see Cartier \cite{C79}, page 145. For item (\ref{Fact2}) see Cartier \cite{C79}, page 140. Item (\ref{Fact3}) can be checked by a direct calculation using equations (3) and (9) on page 145 in Cartier \cite{C79}. The Iwasawa decomposition holds because $K$ is the stabilizer of a special vertex. (See Tits \cite{T79}, 3.3.2.) \EndProof\\

The next proposition states some basic facts about parahoric subgroups. First we introduce some notation. 
Let $\alpha\in \Delta$. Let $I_{\alpha}$ be the parahoric subgroup associated to the facet $\mathscr{F}_{\alpha}$ of $C$ fixed by $s_{\alpha}$. (See Vigneras \cite{V16}, Section 3.7.) Note that $I_{\alpha}$ is contained in the fixator of $\mathscr{F}_{\alpha}$. 

\begin{proposition}\label{Parahoric} $ $
\begin{enumerate}
\item\label{IwahoriFact1} The multiplication map induces an isomorphism 
\begin{equation}
(\prod_{\alpha\in\Phi^{+}}I\cap U_{\alpha})\times(I\cap M)\times(\prod_{\alpha\in\Phi^{-}}I\cap U_{\alpha})\cong I.
\end{equation} 
\noindent(The factors in the product over $\Phi^{\pm}$ may be taken in any order.) 

Let $\alpha\in\Delta$. Then:
\item\label{IwahoriFact2} $I\cap U_{\alpha} = I_{\alpha}\cap U_{\alpha} = K\cap U_{\alpha}$;
\item\label{IwahoriFact3} $I\cap U_{-\alpha}\subsetneq I_{\alpha}\cap U_{-\alpha} = K\cap U_{-\alpha}$;
\item\label{IwahoriFact4} $I\cap \overline{U}\subsetneq I_{\alpha}\cap \overline{U}$;
\item\label{IwahoriFact5} $I_{\alpha}=I\cup Iw_{\alpha}I$.
\end{enumerate}
\end{proposition}

\textbf{Proof:} Item (\ref{IwahoriFact1}) is the Iwahori factorization. For references see Tits \cite{T79} Section 3.1.1 and 3.7, or Vigneras \cite{V16} Section 3.7. Items (\ref{IwahoriFact2}) and (\ref{IwahoriFact3}) follow from  Vigneras \cite{V16} examples 3.10 on page 709 and 3.11 on page 710, and Line (49) page 712. Item (\ref{IwahoriFact4}) follows directly from item (\ref{IwahoriFact3}). Item (\ref{IwahoriFact5}) follows from Vigneras \cite{V16} Proposition 3.26, and Theorem 3.33.\EndProof

\subsection{Function Spaces} 

If $R$ is a $\mathbb{C}$-vector space and $X$ is an $\ell$-space (in the sense of Bernstein-Zelevinski \cite{BZ76}), let $C^{\infty}(X,R)$ be the set of functions $f:X\rightarrow R$ such that $f$ is locally constant and let $C_{c}^{\infty}(X,R)$ be the subset of $C^{\infty}(X,R)$ consisting of compactly supported functions. When $R=\mathbb{C}$ we omit $\mathbb{C}$ from the notation.

Given $f_{1}\in C^{\infty}_{c}(G)$ and $f_{2}\in C^{\infty}_{c}(G,R)$ we define $f_{2}*f_{1}(g)\stackrel{\mathrm{def}}{=}\int_{G}f_{2}(h)f_{1}(h^{-1}g)dg\in C^{\infty}(G,R)$, where $dg$ is the Haar measure on $G$ such that $\mathrm{meas}(I)=1$.

Unless otherwise stated, the group $G$ will act on $(\varrho,C^{\infty}(G,R))$, on the left, via right translation, $(\varrho(h)\cdot f)(g)=f(gh)$. This action induces an action of $C^{\infty}_{c}(G)$ on $C^{\infty}(G,R)$. Let $f_{1}\in C^{\infty}_{c}(G)$ and $f_{2}\in C^{\infty}(G,R)$, then $\varrho(f_{1})\cdot f_{2}=f_{2}*\hat{f}_{1}$, where $\hat{f}_{1}(g)=f_{1}(g^{-1})$.

When $H$ is a subgroup of $G$ and $(\sigma, V)$ is a smooth $H$-representation we define the smooth $G$-representations 
\begin{align*}
\mathrm{Ind}_{H}^{G}(\sigma)=&\{f\in C^{\infty}(G,V)|f(hg)=\sigma(h)f(g)\text{ for all }h\in H\},\\
\mathrm{ind}_{H}^{G}(\sigma)=&\{f\in \mathrm{Ind}_{H}^{G}(\sigma)| \text{ the support of }f\text{ is compact mod }H\},
\end{align*}
where $G$ acts by right translation.

The space $C^{\infty}_{c}(M/M^{\circ})$ is a $\mathbb{C}$-algebra under convolution. For $\lambda\in X_{*}(\textbf{A})$ let $1_{m_{\lambda}M^{\circ}}$ be the characteristic function of the set $m_{\lambda}M^{\circ}$. The set $\{1_{m_{\lambda}M^{\circ}}|\lambda\in X_{*}(\textbf{A})\}$ is a basis for $C^{\infty}_{c}(M/M^{\circ})$.

For $\lambda\in X_{*}(\textbf{A})$ let $e^{\lambda}\in \mathbb{C}[X_{*}(\text{A})]$ be the element associated to $\lambda$. The map $1_{m_{\lambda}M^{\circ}}\mapsto e^{\lambda}$ defines a  $\mathbb{C}$-algebra isomorphism $C^{\infty}_{c}(M/M^{\circ})\rightarrow \mathbb{C}[X_{*}(\textbf{A})]$ since $M/M^{\circ}\cong X_{*}(\textbf{A})$.

\subsection{Satake Transformation}\label{SatTrans} Let $\mathcal{S}:C_{c}^{\infty}(G/K)\rightarrow C_{c}^{\infty}(M/M^{\circ})\cong \mathbb{C}[X_{*}(\textbf{A})]$ be the Satake transform, defined by 
\begin{equation}
\mathcal{S}(f)(m)=\delta_{\overline{P}}(m)^{-1/2}\int_{\overline{U}}f(um)du,
\end{equation}
where the Haar measure of $\overline{U}$ is normalized so that $\mathrm{meas}(\overline{U}\cap K)=1$.

The space $\mathcal{H}_{K}\stackrel{\mathrm{def}}{=} C^{\infty}_{c}(K\backslash G/K)$ has a multiplication defined by $f_{1}*f_{2}(g)=\int_{G}f_{1}(h)f_{2}(h^{-1}g)dh$, where the Haar measure is normalized so that $I$ has measure 1. 

\begin{theorem}[Satake Isomorphism] The map $\mathcal{S}:\mathcal{H}_{K}\rightarrow \mathbb{C}[M/M^{\circ}]^{W}$, defined by
\begin{equation}
\mathcal{S}(f)(m)=\delta_{\overline{P}}(m)^{-1/2}\int_{\overline{U}}f(um)du,
\end{equation}
is an isomorphism of $\mathbb{C}$-algebras.
\end{theorem}
\noindent(For details see Cartier \cite{C79} Theorem 4.1, page 147.) 

A direct calculation shows that $\mathcal{S}:C_{c}^{\infty}(G/K)\rightarrow C_{c}^{\infty}(\overline{U}\backslash G/K)\cong \mathbb{C}[M/M^{\circ}]$ is an $\mathcal{H}_{K}$-module homomorphism where $\mathcal{H}_{K}$ acts on $\mathbb{C}[M/M^{\circ}]$ through the Satake isomorphism.

\subsection{Universal Principal Series}

Our presentation of the Casselman-Shalika formula will utilize the universal principal series, which we now introduce. For additional details see the survey article of Haines-Kottwitz-Prasad \cite{HKP10}. For the remainder of the paper let $R=\mathbb{C}[M/M^{\circ}]\cong\mathbb{C}[X_{*}(\textbf{A})]$ (where the isomorphism is induced via $\nu$) and define $\chi_{\mathrm{univ}}:M\rightarrow R^{\times}$ to be the tautological character. The universal principal series is defined to be 
\begin{equation}
i_{G,\overline{P}}(\chi_{\mathrm{univ}}^{-1})=\mathrm{Ind}_{\overline{P}}^{G}(\delta_{\overline{P}}^{1/2}\otimes\chi_{\mathrm{univ}}^{-1})=\{f\in C^{\infty}(G, R)|f(mug)=(\delta_{\overline{P}}^{1/2}\cdot\chi_{\mathrm{univ}}^{-1})(m)f(g)\}.
\end{equation}
The space $i_{G,\overline{P}}(\chi_{\mathrm{univ}}^{-1})$ is an $(R,G)$-bimodule, where $r\in R$ acts on $f\in i_{G,\overline{P}}(\chi_{\mathrm{univ}}^{-1})$ via $(r\cdot f)(g)=r\cdot(f(g))$; $G$ acts via right translation.

Let $v_{0}\in i_{G,\overline{P}}(\chi_{\mathrm{univ}}^{-1})^{K}$ be the spherical vector normalized such that $v_{0}(1)=1$.

It will be convenient to extend the scalars of the universal principal series. Let $S\supseteq R$ be a commutative $\mathbb{C}$-algebra that is an integral domain. Write $i_{S}:R\rightarrow S$ for the inclusion of $R$ into $S$.

\begin{lemma}\label{ExtScalePSeries}
The map $\varpi:S\otimes_{R} i_{G,\overline{P}}(\chi_{\mathrm{univ}}^{-1})\rightarrow  i_{G,\overline{P}}(i_{S}\circ\chi_{\mathrm{univ}}^{-1})$ defined by $s\otimes f(g)\mapsto si_{S}(f(g))$ is an isomorphism of $(S,G)$-bimodules.
\end{lemma}

\textbf{Proof:} This follows from smoothness and because $S$ is an integral domain.\EndProof

\subsection{Relating Universal Principal Series}\label{RUPS}

In this subsection we describe two situations in which we can relate the universal principal series on two groups. First we discuss the case where $\textbf{G}^{\prime\prime}=\textbf{G}\times_{\mathrm{Spec}(F)}\textbf{T}$, where $\textbf{G}$ is a  connected unramified group and $\textbf{T}$ is a torus. The group scheme $\textbf{G}^{\prime\prime}$ has a maximal split torus $\textbf{A}^{\prime\prime}=\textbf{A}\times_{\mathrm{Spec}(F)} \textbf{T}$ contained in the minimal parabolic subgroup $\textbf{P}^{\prime\prime}=\textbf{P}\times_{\mathrm{Spec}(F)} \textbf{T}$ with Levi subgroup $\textbf{M}^{\prime\prime}=\textbf{M}\times_{\mathrm{Spec}(F)} \textbf{T}$.

Let $R^{\prime\prime}=\mathbb{C}[M^{\prime\prime}/(M^{\prime\prime})^{\circ}]$ and let $\xi:T\rightarrow \mathbb{C}[T/T^{\circ}]^{\times}$ be the tautological character of $T$. Note that the tautological character $\chi_{\mathrm{univ}}^{\prime\prime}:M^{\prime\prime}\rightarrow R^{\prime\prime}\cong R\otimes_{\mathbb{C}} \mathbb{C}[T/T^{\circ}]$ of the Levi subgroup $M^{\prime\prime}$ is given by $\chi_{\mathrm{univ}}^{\prime\prime}((m,t))=\chi_{\mathrm{univ}}(m)\otimes \xi(t)$.

\begin{lemma}\label{ProdPSeries}
 With the notation above, the map $\tau:R^{\prime\prime}\otimes_{R}i_{G,\overline{P}}(\chi_{\mathrm{univ}}^{-1})\rightarrow i_{G^{\prime\prime},\overline{P}^{\prime\prime}}((\chi_{\mathrm{univ}}^{\prime\prime})^{-1})$ defined by $\chi\otimes f\mapsto ((g,t)\mapsto \chi\xi^{-1}(t)f(g))$ is an isomorphism of $(R^{\prime\prime},G^{\prime\prime})$-bimodules.
\end{lemma}

\textbf{Proof:} This is a consequence of smoothness.\EndProof\\

\textbf{Remark:} Recall that the right regular action of $T$ on $\mathbb{C}[T/T^{\circ}]$ when its elements are viewed as functions corresponds to multiplication by $\xi^{-1}$ when $\mathbb{C}[T/T^{\circ}]$ is viewed as a group algebra.\\

Second, we consider the case of two connected unramified groups $\textbf{G}^{\prime\prime}$ and $\textbf{G}^{\prime}$ and a map $\pi:G^{\prime}\rightarrow G^{\prime\prime}$ that is an open algebraic group homomorphism with finite central kernel such that $P^{\prime\prime}\mathrm{Im}(\pi)=G^{\prime\prime}$. In this case, we take $P^{\prime}=\pi^{-1}(P^{\prime\prime})$ as our minimal parabolic subgroup of $G^{\prime}$. This implies that $\pi(M^{\prime})\subseteq M^{\prime\prime}$, $\pi(A^{\prime})\subseteq A^{\prime\prime}$, and $\pi(U^{\prime})\subseteq U^{\prime\prime}$.

Let $R^{\prime\prime}=\mathbb{C}[X_{*}(\textbf{A}^{\prime\prime})]$ and $R^{\prime}=\mathbb{C}[X_{*}(\textbf{A}^{\prime})]$. Then we have the following maps:
\begin{align}
\pi_{*}:&X_{*}(\textbf{A}^{\prime})\hookrightarrow X_{*}(\textbf{A}^{\prime\prime}) &(\ker\pi\text{ finite implies injectivity}),\label{cocharinj}\\
\pi^{*}:&i_{G^{\prime\prime},\overline{P}^{\prime\prime}}((\chi_{\mathrm{univ}}^{\prime\prime})^{-1})\rightarrow\, i_{G^{\prime},\overline{P}^{\prime}}(\pi_{*}\circ(\chi_{\mathrm{univ}}^{\prime})^{-1}),
\end{align}
where $\pi^{*}(f)(g^{\prime})=f(\pi(g^{\prime}))$. The map $\pi^{*}$ is well-defined because 
\begin{equation*}
\pi^{*}(f)(m^{\prime}g^{\prime})=(\delta_{\overline{P}^{\prime\prime}}^{1/2}\chi_{\mathrm{univ}}^{-1})(\pi(m^{\prime}))\pi^{*}(f)(g^{\prime}),
\end{equation*}
$(\chi_{\mathrm{univ}}^{\prime\prime})^{-1}\circ \pi_{*} = \pi_{*} \circ (\chi_{\mathrm{univ}}^{\prime})^{-1}$, and $\delta_{\overline{P}^{\prime\prime}}(\pi(m^{\prime}))=\delta_{\overline{P}^{\prime}}(m^{\prime})$. Furthermore, if $g^{\prime\prime}=\pi(g^{\prime})$, then  
\begin{equation}
\pi^{*}(\varrho(g^{\prime\prime})f)=\varrho(g^{\prime})\pi^{*}(f).
\end{equation}

\begin{lemma}\label{GroupTransfer1}
The map $\pi^{*}:i_{G^{\prime\prime},\overline{P}^{\prime\prime}}((\chi_{\mathrm{univ}}^{\prime\prime})^{-1})\rightarrow\, i_{G^{\prime},\overline{P}^{\prime}}(\pi_{*}\circ(\chi_{\mathrm{univ}}^{\prime})^{-1})$ is an isomorphism of $(R^{\prime\prime},G^{\prime})$-bimodules.
\end{lemma} 

\textbf{Proof:} We show that the map is a bijection. The rest follows directly from definitions. Since $P^{\prime\prime}\mathrm{Im}(\pi)=G^{\prime\prime}$, the map $\pi^{*}$ is injective.

Now we show that $\pi^{*}$ is surjective. Let $x^{\prime}\in G^{\prime}$ and $K^{\prime}\subset G^{\prime}$ be an open compact subgroup such that $K^{\prime}\cap (x^{\prime})^{-1}P^{\prime}x^{\prime}\subseteq (x^{\prime})^{-1}(M^{\prime})^{\circ}U^{\prime}x^{\prime}$. Define $\phi_{x^{\prime},K^{\prime}}^{\prime}$ be the function supported on $P^{\prime}x^{\prime}K^{\prime}$ such that $\phi(m^{\prime}u^{\prime}x^{\prime}k^{\prime})=\pi_{*}\circ(\chi_{\mathrm{univ}}^{\prime})^{-1}(m^{\prime})$, where $m^{\prime}\in M^{\prime}$, $u^{\prime}\in \overline{U}$, $k^{\prime}\in K^{\prime}$. The module $i_{G^{\prime},\overline{P}^{\prime}}(\pi_{*}\circ(\chi_{\mathrm{univ}}^{\prime})^{-1})$ is the $R^{\prime\prime}$-span of the functions $\phi_{x^{\prime},K^{\prime}}^{\prime}$.

Similarly, we can define $\phi_{\pi(x^{\prime}),\pi(K^{\prime})}\in i_{G^{\prime\prime},\overline{P}^{\prime\prime}}((\chi_{\mathrm{univ}}^{\prime\prime})^{-1})$. (Since $\pi$ is an open map, $\pi(K^{\prime})$ is a compact open subgroup of $G^{\prime\prime}$.) By construction $\pi^{*}(\phi_{\pi(x^{\prime}),\pi(K^{\prime})})=\phi_{x^{\prime},K^{\prime}}$. Thus $\pi^{*}$ is surjective.\EndProof

We Combine Lemmas \ref{ExtScalePSeries} (with $S=R^{\prime\prime}$) and \ref{GroupTransfer1} to get the following. 

\begin{lemma}\label{GroupTransfer2}
The map $\varpi^{-1}\circ\pi^{*}:i_{G^{\prime\prime},\overline{P}^{\prime\prime}}((\chi_{\mathrm{univ}}^{\prime\prime})^{-1})\rightarrow R^{\prime\prime}\otimes_{R^{\prime}} i_{G^{\prime},\overline{P}^{\prime}}((\chi_{\mathrm{univ}}^{\prime})^{-1})$ is an isomorphism of $(R^{\prime\prime},G^{\prime})$-bimodules.
\end{lemma}

We use the above results in the following setting. Let $\textbf{G}^{\prime}$ be a connected unramified group. Define the torus $\textbf{T}=\textbf{G}^{\prime}/\mathcal{D}\textbf{G}^{\prime}$ and the connected semisimple adjoint group $\textbf{G}=\textbf{G}^{\prime}/Z(\textbf{G}^{\prime})$ and let $\boldsymbol{\pi}:\textbf{G}^{\prime}\rightarrow \textbf{G}^{\prime\prime}=\textbf{G}\times_{\mathrm{Spec}(F)} \textbf{T}$ be the natural map induced by the quotient maps.  The map $\boldsymbol{\pi}$ induces a map on $F$-points $\pi:G^{\prime}\rightarrow G\times T$ such that $\pi$ is an open algebraic group homomorphism with finite central kernel such that $(P\times T)\mathrm{Im}(\pi)=G\times T$. In this case, we will take $S=R^{\prime\prime}$.

\subsection{Whittaker Function}\label{WF}

Let $\Psi:\overline{U}\rightarrow \mathbb{C}^{\times}$ be a smooth character. The character $\Psi$ factors through $\overline{U}/[\overline{U},\overline{U}]\cong \prod_{\alpha\in -\Delta}U_{\alpha}/U_{2\alpha}$, where $U_{2\alpha}=\{1\}$ if $2\alpha\notin \Phi$. Thus to define $\Psi$ it suffices to define smooth characters $\Psi_{\alpha}:U_{\alpha}\rightarrow \mathbb{C}^{\times}$ for all $\alpha\in -\Delta$ and set $\Psi=\prod_{\alpha\in -\Delta}\Psi_{\alpha}$. We say that $\Psi$ is nondegenerate if $\Psi|_{U_{\alpha}}$ is not trivial for all $\alpha\in -\Delta$. We write $\overline{\Psi}$ for the complex conjugate of $\Psi$. 

We are interested in $R$-valued Whittaker functionals, but it will be convenient to allow for extension of scalars. Again we let $S\supseteq R$ be a $\mathbb{C}$-algebra that is an integral domain with an inclusion $i_{S}:R\rightarrow S$. For any abelian group $J$, define $J_{\Psi}$ to be the abelian group $J$ with an action of $\overline{U}$ given by $u\cdot j=\Psi(u)j$. Let 
\begin{equation*}
i_{G,\overline{P}}(i_{S}\circ\chi_{\mathrm{univ}}^{-1})_{(\overline{U},\Psi)}=i_{G,\overline{P}}(i_{S}\circ\chi_{\mathrm{univ}}^{-1})/\mathrm{span}_{S}(u\cdot f-\Psi(u)f|u\in \overline{U},f\in i_{G,\overline{P}}(i_{S}\circ\chi_{\mathrm{univ}}^{-1}))
\end{equation*}
be the $\Psi$-twisted Jacquet module of $i_{G,\overline{P}}(i_{S}\circ\chi_{\mathrm{univ}}^{-1})$. An $S$-valued Whittaker functional is an element $\mathfrak{W}\in \mathrm{Hom}_{(S,\overline{U})}(i_{G,\overline{P}}(i_{S}\circ\chi_{\mathrm{univ}}^{-1}),S_{\Psi})$. Note that for any $S$-module $S^{\prime}$ we have 
\begin{equation}\label{WhittFunctional}
\mathrm{Hom}_{(S,\overline{U})}(i_{G,\overline{P}}(i_{S}\circ\chi_{\mathrm{univ}}^{-1}),S^{\prime}_{\Psi})\cong \mathrm{Hom}_{S}(i_{G,\overline{P}}(i_{S}\circ\chi_{\mathrm{univ}}^{-1})_{(\overline{U},\Psi)},S^{\prime})
\end{equation}
as $S^{\prime}$-modules. 

\begin{lemma}\label{WhittFunctExtScale}
The natural map
\begin{equation}
\mathrm{Hom}_{(S,\overline{U})}(S\otimes_{R}i_{G,\overline{P}}(\chi_{\mathrm{univ}}^{-1}),S_{\Psi})\rightarrow \mathrm{Hom}_{(R,\overline{U})}(i_{G,\overline{P}}(\chi_{\mathrm{univ}}^{-1}),S_{\Psi})
\end{equation}
defined by $\phi\mapsto \phi|_{1\otimes i_{G,\overline{P}}(\chi_{\mathrm{univ}}^{-1})}$ is an isomorphism of $S$-modules.
\end{lemma}

\textbf{Proof:} The inverse map is defined by $\phi^{\prime}\mapsto (s\otimes f\mapsto s\phi^{\prime}(f))$.\EndProof\\

The $S$-valued spherical Whittaker function associated to $\mathfrak{W}$ is defined to be $\mathcal{W}(g)=\mathfrak{W}(gv_{0})$. The Iwasawa decomposition $G=\overline{U}MK$ shows that $\mathcal{W}$ is determined by its values on $M/M^{\circ}\cong A/A^{\circ}$. Note that $\mathcal{W}\in \mathrm{Ind}_{\overline{U}}^{G}(S_{\Psi})^{K}$. The Casselman-Shalika formula is a formula for the function $\mathcal{W}$ evaluated at points in $A/A^{\circ}$.

\subsection{Dual Group}\label{DualGroup}
In this subsection we introduce a complex group whose characters will appear in the Casselman-Shalika formula. Let $\mathscr{X}= X_{*}(\textbf{A})+\Lambda^{\vee}\subseteq X_{*}(\textbf{A})\otimes_{\mathbb{Z}} \mathbb{Q}$ and $\mathscr{Y}=\mathrm{Hom}_{\mathbb{Z}}(\mathscr{X},\mathbb{Z})$. Let $\Phi^{\mathrm{nd}}$ be the set of non-divisible roots in $\Phi$ and let $(\Phi^{\mathrm{nd}})^{\vee}$ be the set of coroots of $\Phi^{\mathrm{nd}}$.

\begin{proposition}
The quadruple $(\mathscr{X},(\Phi^{\mathrm{nd}})^{\vee},\mathscr{Y},\Phi^{\mathrm{nd}})$ defines a root datum.
\end{proposition}

\textbf{Proof:} The only point that may not be clear is that $\Phi^{\mathrm{nd}}\subset \mathscr{Y}$. Let $\alpha\in \Phi^{\mathrm{nd}}$. Then by definition $\alpha\in \mathrm{Hom}_{\mathbb{Z}}(X_{*}(\textbf{A}),\mathbb{Z})$. Thus $\alpha$ extends uniquely to an element of $\mathrm{Hom}_{\mathbb{Q}}(X_{*}(\textbf{A})\otimes_{\mathbb{Z}}\mathbb{Q},\mathbb{Q})$. Since $\mathscr{X}\subset X_{*}(\textbf{A})\otimes_{\mathbb{Z}}\mathbb{Q}$ and $\alpha\in \mathrm{Hom}_{\mathbb{Z}}(\Lambda^{\vee},\mathbb{Z})$ it follows that $\alpha\in \mathrm{Hom}_{\mathbb{Z}}(\mathscr{X},\mathbb{Z})$. \EndProof\\

\noindent Let ${^{L}\textbf{G}}^{\dagger}$ be the connected complex reductive group with root datum $(\mathscr{X},(\Phi^{\mathrm{nd}})^{\vee},\mathscr{Y},\Phi^{\mathrm{nd}})$. Let $\rho^{\vee}$ be half the sum of the positive roots (with respect to $P$) in $(\Phi^{\mathrm{nd}})^{\vee}$. Note that $\rho^{\vee}\in \mathscr{X}$.

\textbf{Examples of ${^{L}\textbf{G}}^{\dagger}$:} 
\begin{enumerate}
\item If $G$ is a semisimple split group, then $\mathscr{X}$ is the coweight lattice of $G$ and $\Phi^{\mathrm{nd}}=\Phi$. Thus ${^{L}\textbf{G}}^{\dagger}$ is the simply connected cover of the Langlands dual group of $G$.\\
\item Let $E/F$ be an unramified quadratic extension. Let $\textbf{G}=\mathrm{SU}(2n+1,E/F)$ be the unramified quasi-split special unitary group over $F$ associated to the Hermitian form 
$\left(\begin{smallmatrix}
&&1\\
&\iddots&\\
1&&\end{smallmatrix}\right)$ on the vector space $E^{2n+1}$. In this case, ${^{L}\textbf{G}}^{\dagger}=\mathrm{Sp}(2n,\mathbb{C})$. 
\end{enumerate}

Let ${^{L}\textbf{A}}=X^{*}(\textbf{A})\otimes \mathbb{C}^{\times}$ and let ${^{L}\textbf{A}}^{\dagger}$ be a maximal torus of $^{L}G^{\dagger}$. Since ${^{L}\textbf{A}}^{\dagger}\cong \mathscr{Y}\otimes \mathbb{C}^{\times}$, the restriction map $\mathscr{Y}\rightarrow X^{*}(\textbf{A})$ induces an isogeny ${^{L}\textbf{A}}^{\dagger}\rightarrow {^{L}\textbf{A}}$ of tori. In addition we have the inclusion $X_{*}(\textbf{A})\hookrightarrow\mathscr{X}$.

For $H$ a complex reductive group, let $\mathrm{Rep}(H)$ denote the $\mathbb{C}$-algebra of characters of finite dimensional algebraic representations of $H$. 

\begin{lemma} The inclusion $X^{*}({^{L}\textbf{\textbf{A}}})=X_{*}(\textbf{\textbf{A}})\stackrel{inc}{\hookrightarrow} \mathscr{X}$ induces an inclusion
$\mathbb{C}[{^{L}\textbf{\textbf{A}}}/W]\hookrightarrow\mathbb{C}[{^{L}\textbf{A}}^{\dagger}/W]\cong \mathrm{Rep}({^{L}\textbf{G}}^{\dagger})$.
\end{lemma}

\textbf{Proof:} ${^{L}\textbf{G}}^{\dagger}$ is a connected reductive group over $\mathbb{C}$, thus $\mathrm{Rep}({^{L}\textbf{G}}^{\dagger})\cong \mathbb{C}[{^{L}\textbf{\textbf{A}}}^{\dagger}/W]$.\EndProof\\

For any $\lambda\in X_{*}({^{L}\textbf{A}})^{+}$ define $A_{\lambda}$ to be the element of  $\mathcal{H}_{K}$ that corresponds to $\mathrm{ch}V_{inc(\lambda)}$ under the inclusion $\mathcal{H}_{K}\stackrel{\mathcal{S}}{\rightarrow}\mathbb{C}[{^{L}\textbf{A}}/W]\hookrightarrow \mathrm{Rep}({^{L}\textbf{G}}^{\dagger})$. Later we will suppress $inc$ in our notation.

The following coefficients play an important role in our proof of the Casselman-Shalika formula. Given $\lambda\in \mathscr{X}^{+}$ and $\mu,\eta\in \mathscr{X}^{++}$ define $c_{\lambda,\mu}^{\eta}$ so that 
\begin{equation}
\mathrm{ch}V_{\lambda}\cdot \mathrm{ch}V_{\mu-\rho^{\vee}}=\sum_{\eta}c_{\lambda,\mu}^{\eta}\mathrm{ch}V_{\eta-\rho^{\vee}}.
\end{equation} 
Recall that if $c_{\lambda,\mu}^{\eta}\neq 0$, then $\lambda+\mu-\eta\in \mathrm{span}_{\mathbb{Z}}((\Phi^{\mathrm{nd}})^{\vee})$.


\section{Iwahori-Hecke Algebra}\label{IwahoriHecke}

In this section we recall the Bernstein presentation of the Iwahori-Hecke algebra $\mathcal{H}\stackrel{\mathrm{def}}{=}\mathcal{H}_{I}=C_{c}^{\infty}(I \backslash G/I)$ following the presentation of Rostami \cite{R15}. We begin with some preliminaries. The space $\mathcal{H}$ has a multiplication defined by $f_{1}*f_{2}(g)=\int_{G}f_{1}(h)f_{2}(h^{-1}g)dh$ and the Haar measure is normalized so that $I$ has measure 1. Let $1_{X}$ be the characteristic function of the set $X\subset G$, and for $g\in G$ let $\mathcal{T}_{g}=1_{IgI}\in\mathcal{H}$. It is convenient to normalize $\mathcal{T}_{g}$, so let $\textbf{q}(g)\stackrel{\mathrm{def}}{=}[IgI:I]$ and define $\overline{\mathcal{T}}_{g}\stackrel{\mathrm{def}}{=}\textbf{q}(g)^{-1/2}\mathcal{T}_{g}$.

Next we introduce the elements of Bernstein's commutative subalgebra. For $\lambda\in X_{*}(\textbf{A})$ we choose $\lambda_{1},\lambda_{2}\in X_{*}(\textbf{A})^{+}$ such that $\lambda=\lambda_{1}-\lambda_{2}$ and define
\begin{equation}
\theta_{\lambda}\stackrel{\mathrm{def}}{=}\overline{\mathcal{T}}_{m_{\lambda_{1}}}*\overline{\mathcal{T}}_{m_{\lambda_{2}}}^{-1}.
\end{equation}
The definition of $\theta_{\lambda}$ depends neither on our choice of the $\lambda_{j}$ nor the splitting $s$ of the exact sequence on line (\ref{ExactSq1}). The invertibility of $\mathcal{T}_{m_{\lambda}}$ follows from the first part of Theorem \ref{BP}; the elements $\theta_{\lambda}$ do not appear until the second part of that theorem. The elements $\theta_{\lambda}$ generate a commutative $\mathbb{C}$-algebra $\mathcal{A}$, which is isomorphic to $\mathbb{C}[X_{*}(\textbf{A})]$.

Let $\widetilde{W}=N/M^{\circ}$ be the extended affine Weyl group. It acts on $\mathscr{A}$ (Tits \cite{T79}, Section 1.2) and the stablizer of $x_{0}$ is isomorphic to $W$, the Weyl group of the relative root system $\Phi$. For each $w\in W$ we will choose a representative $n_{w}\in N\cap\mathrm{Stab}_{G}(x_{0})=N\cap K$. We write $\mathcal{T}_{w}$ in place of $\mathcal{T}_{n_{w}}$. Note that $\mathcal{T}_{w}$ is independent of the choice of $n_{w}$ because $M^{\circ}\subseteq I$.

If $H$ is the group of $F$-points of a connected reductive group scheme, then we let $\Omega_{H}$ be the image of the Kottwitz homomorphism (\cite{K97}, Sections 7.1-7.4). Now we can state the Iwahori-Matsumoto presentation and the Bernstein presentation for $\mathcal{H}$. (Note that because $G$ is quasi-split, $M$ is a torus and $\Omega_{M}=X_{*}(\textbf{A})$.)

\begin{theorem}[Rostami \cite{R15}; Vign\'{e}ras \cite{V16}]\label{BP} The product of any two basis elements $\mathcal{T}_{w}$ is determined by the relations 
\begin{itemize}
\item $\mathcal{T}_{w}*\mathcal{T}_{w^{\prime}}=\mathcal{T}_{ww^{\prime}}$, for any $w,w^{\prime}\in \widetilde{W}$, such that  $\ell(w w^{\prime})=\ell(w)+\ell(w^{\prime})$;\label{LengthMult}
\item $\mathcal{T}_{s}*\mathcal{T}_{s} = (q(s)-1)\mathcal{T}_{s}+q(s)$, for any simple affine reflection $s\in \Delta_{\mathrm{aff}}$;\label{QuadRel}
\item $\mathcal{T}_{w}*\mathcal{T}_{\tau}=\mathcal{T}_{w\tau}=\mathcal{T}_{\tau}*\mathcal{T}_{\tau^{-1}w\tau}$, for any $w\in \widetilde{W}, \tau\in \Omega_{G}$. \label{TwistMult}
\end{itemize}

The set $\{\mathcal{T}_{w}*\theta_{\lambda}|\lambda\in \Omega_{M},w\in W\}$ is a $\mathbb{C}$-linear basis for the Iwahori-Hecke algebra $\mathcal{H}$, and the product of two basis elements is determined by the relations:
\begin{itemize}
\item[($Add$)] $\theta_{\lambda}*\theta_{\mu}=\theta_{\lambda+\mu}$, for any $\lambda,\mu\in \Omega_{M}$.\label{Add}
\end{itemize}
\noindent If $\alpha\in \Delta$ and $\lambda\in \Omega_{M}$, then
\begin{itemize}
\item[($BR$)] $\mathcal{T}_{s_{\alpha}}*\theta_{\lambda}=$
$\begin{cases}
\theta_{s_{\alpha}(\lambda)}*\mathcal{T}_{s_{\alpha}}+\sum_{j=0}^{\langle\alpha,\lambda\rangle-1}\textbf{q}_{j}(s_{\alpha})\theta_{\lambda-j\alpha^{\vee}}& \text{ if }\langle \alpha,\lambda\rangle> 0;\\
\theta_{s_{\alpha}(\lambda)}*\mathcal{T}_{s_{\alpha}}&\text{ if }\langle \alpha,\lambda\rangle= 0;\\
\theta_{s_{\alpha}(\lambda)}*\mathcal{T}_{s_{\alpha}}-\sum_{j=0}^{\langle\alpha,s_{\alpha}(\lambda)\rangle-1}\textbf{q}_{j}(s_{\alpha})\theta_{s_{\alpha}(\lambda)-j\alpha^{\vee}}& \text{ if }\langle \alpha,\lambda\rangle<0.
\end{cases}$\label{BR}
\end{itemize}
(We will not need the definition of $\textbf{q}_{j}(s)$, which can be found in Rostami \cite{R15}, Section 5.4.)
\end{theorem}
$ $\\
\noindent\textbf{Remark:} The Bernstein relations (BR) above include a minor correction to what appears in Rostami \cite{R15}. Namely Rostami's formula only holds in the case $\langle\alpha,\lambda\rangle\geq 0$. Fortunately, the case $\langle\alpha,\lambda\rangle< 0$ can be deduced directly from the case $\langle\alpha,-\lambda\rangle> 0$.\\

One consequence of Theorem \ref{BP} is that the intermediate algebra $\mathcal{H}_{IK}\stackrel{\mathrm{def}}{=}C^{\infty}_{c}(I\backslash G/K)$ has a nice basis. Since $K=\cup_{w\in W}IwI$, the characteristic function of $K$ can be expressed as $1_{K} = \sum_{w\in W}\mathcal{T}_{w}$. This identity, Theorem \ref{BP}, and the fact that $f(g)\mapsto f(g^{-1})$ defines an anti-isomorphism of $\mathcal{H}$ yield the following corollary.
\begin{corollary}
The Hecke algebra $\mathcal{H}_{IK}$ has a basis consisting of the elements $\theta_{\lambda}^{K}\stackrel{\mathrm{def}}{=}\theta_{\lambda}*1_{K}$, where $\lambda\in X_{*}(\textbf{A})$.
\end{corollary}

We end this section with an identity that is crucial for the proof of Proposition \ref{TwistedKer}. It follows directly from the Bernstein relations.

\begin{corollary}\label{CommId} Let $\alpha\in \Delta$ and let $s$ be the simple reflection associated to $\alpha$. Let $\lambda\in X_{*}(\textbf{A})$. Then
\begin{equation}
\mathcal{T}_{s}(\theta_{\lambda}+\theta_{s(\lambda)})=(\theta_{\lambda}+\theta_{s(\lambda)})\mathcal{T}_{s}.
\end{equation}
\end{corollary}

\section{Savin's Isomorphism}\label{SavinsIsoSec}

In this section we prove Savin's Isomorphism (\cite{S12}, Theorem 1) for connected unramified groups, which states that the Satake transform defines an $\mathcal{H}_{K}$-module isomorphism from $\mathcal{H}_{IK}$ to $C^{\infty}_{c}(M/M^{\circ})\cong \mathbb{C}[X_{*}(\textbf{A})]$. Using the results of Section \ref{IwahoriHecke} this amounts to a minor modification of Savin \cite{S12}. We begin with some notation and then a few basic results.

Given an equation
\begin{equation}\label{Test}
A=B,
\end{equation}
let $\mathrm{RHS}(\ref{Test})=B$ ($\mathrm{LHS}(\ref{Test})=A$) be the right(left)-hand-side of equation (\ref{Test}).

\begin{lemma}\label{EquiVar}
Let $(\pi,V)$ be a smooth $G$-module and $(\pi^{\prime},V^{\prime})$ a smooth $\overline{P}$-module with the trivial action of $\overline{U}$. Let $\mathcal{S}:V\rightarrow V^{\prime}$ be a map such that $\mathcal{S}(\pi(p)v)=\delta_{\overline{P}}^{1/2}(p)\pi^{\prime}(p)\mathcal{S}(v)$ for every $p\in \overline{P}$. Then, for every $\lambda\in X_{*}(\textbf{A})$ and $v\in V^{I}$,
\begin{equation}
\mathcal{S}(\pi(\theta_{\lambda})v)=\pi^{\prime}(m_{\lambda})\mathcal{S}(v).
\end{equation}
\end{lemma}

\textbf{Proof:} To begin, we prove this for $\lambda\in X_{*}(\textbf{A})^{+}$, so $\theta_{\lambda}=\delta_{\overline{P}}^{1/2}(m_{\lambda})1_{Im_{\lambda}I}$. By definition

\begin{equation}\label{EquiVarEqn1}
\mathcal{S}(\theta_{\lambda}v) = \delta_{\overline{P}}^{1/2}(m_{\lambda})\mathcal{S}(\int_{Im_{\lambda}I}gvdg).
\end{equation}
Since $v\in V^{I}$ and $\mathrm{meas}(I)=1$ we have
\begin{equation}\label{EquiVarEqn2}
\mathrm{RHS}(\ref{EquiVarEqn1})=\delta_{\overline{P}}^{1/2}(m_{\lambda})\mathcal{S}(\sum_{\gamma\in I/I\cap m_{\lambda}Im_{\lambda}^{-1}}\gamma m_{\lambda}v) .
\end{equation}
By the Iwahori factorization (Proposition \ref{Parahoric}) and the fact that $m_{\lambda}$ contracts $\overline{U}$ (Proposition \ref{ConjContract}) we can assume that each $\gamma\in I/I\cap m_{\lambda}Im_{\lambda}^{-1}$ is represented by an element in $I\cap \overline{U}$. Thus
\begin{equation}\label{EquiVarEqn3}
\mathrm{RHS}(\ref{EquiVarEqn2})=\delta_{\overline{P}}(m_{\lambda})\sum_{\gamma\in I/I\cap m_{\lambda}Im_{\lambda}^{-1}} m_{\lambda}\mathcal{S}(v)=\delta_{\overline{P}}(m_{\lambda})[Im_{\lambda}I:I]m_{\lambda}\mathcal{S}(v).
\end{equation}
Since $[Im_{\lambda}I:I]=\delta_{\overline{P}}^{-1}(m_{\lambda})$ (Proposition \ref{FactsProp1}), we get
\begin{equation}\label{EquiVarEqn4}
\mathrm{RHS}(\ref{EquiVarEqn3})=m_{\lambda}\mathcal{S}(v).
\end{equation}

The general result follows from this special case and is left to the reader. \EndProof\\

\textbf{Remark:} Lemma \ref{EquiVar} is independent of the choice of the splitting $s$ of sequence (\ref{ExactSq1}) because $I\cap M=M^{\circ}$ (Cartier \cite{C79}, page 140).\\

\begin{theorem}[Savin's Isomorphism]\label{SavinIso}
Let $\lambda\in X_{*}(\textbf{A})$. The Satake transform \newline$\mathcal{S}:C^{\infty}_{c}(G/K)\rightarrow C^{\infty}_{c}(M/M^{\circ})\cong \mathbb{C}[X_{*}(\textbf{A})]$ sends the element $\theta_{\lambda}^{K}$ to $e^{\lambda}$. Hence $\mathcal{S}$ induces an isomorphism of left $\mathcal{H}_{K}\cong \mathbb{C}[X_{*}(\textbf{A})]^{W}$-modules
\begin{equation}
\mathcal{H}_{IK}\cong \mathbb{C}[X_{*}(\textbf{A})].
\end{equation}
\end{theorem}

\textbf{Proof:} First note that $\mathcal{S}(1_{K})=1_{M^{\circ}}=e^{0}$. To compute $\mathcal{S}(\theta_{\lambda}^{K})$ we apply Lemma \ref{EquiVar}. In particular, $(\pi,V)$ is the $G$-representation where $V=C_{c}^{\infty}(G/K)$ and $G$ acts by left translation; $(\pi^{\prime},V^{\prime})$ is the $P$-representation where $V^{\prime}=C_{c}(M/M^{\circ})\cong X_{*}(\textbf{A})$ and $P$ acts through the quotient $M\cong P/U$ by left translation. Lemma \ref{EquiVar} states that for any $f\in V^{I}=C^{\infty}_{c}(I\backslash G/K)=\mathcal{H}_{IK}$ we have 

\begin{equation}
\mathcal{S}(\theta_{\lambda}*f)=\pi^{\prime}(m_{\lambda})\mathcal{S}(f).
\end{equation}
In particular, for $f=1_{K}$ we have 

\begin{equation}
\mathcal{S}(\theta_{\lambda})=\mathcal{S}(\theta_{\lambda}*1_{K})=\pi^{\prime}(m_{\lambda})\mathcal{S}(1_{K})=\pi^{\prime}(m_{\lambda})1_{M^{\circ}}=e^{\lambda}.
\end{equation}
\EndProof


\section{Twisted Satake Transform}\label{TwistedSatakeSec}

In this section we study the structure of $\mathrm{ind}_{\overline{U}}^{G}(\overline{\Psi})^{K}$ as an $\mathcal{H}_{K}$-module for a non-degenerate character $\Psi$ of conductor $\frak{p}$. This is accomplished in Proposition \ref{indIso}.

When we say that $\Psi$ has `conductor $\frak{p}$', we mean that for any $\alpha\in \Delta$ the character $\Psi|_{U_{-\alpha}}=\Psi_{\alpha}$ is nontrivial on $U_{-\alpha}\cap I_{\alpha}=U_{-\alpha}\cap K$ and trivial on $U_{-\alpha}\cap I$. One can show that such a character exists by reducing to the case where $G$ is a simply-connected semi-simple unramified group of rank one. (i.e. $G$ is isomorphic to $\mathrm{SL}(2,L)$, where $L$ is an unramified extension of $F$, or $\mathrm{SU}(3,E/F)$, where $E$ is an unramified quadratic extension of $F$.)

We begin our study of $\mathrm{ind}_{\overline{U}}^{G}(\overline{\Psi})^{K}$ by constructing a geometric basis. 
\begin{lemma}\label{CompIndSup}
Suppose that $f\in \mathrm{ind}_{\overline{U}}^{G}(\overline{\Psi})^{K}$. Then for any $\mu\in X_{*}(\textbf{A})\setminus X_{*}(\textbf{A})^{++}$ we have $f(m_{\mu})=0$.
\end{lemma}
\textbf{Proof:} The Iwasawa decomposition shows that it suffices to consider $m_{\mu}$ such that $\mu\in X_{*}(\textbf{A})$. We will show that if $\mu\in X_{*}(\textbf{A})$ is not strictly dominant, then $f(m_{\mu})=0$.

In this case, there exists $\alpha\in -\Delta$ such that $m_{\mu}(U_{\alpha}\cap K)m_{\mu}^{-1}\supseteq (U_{\alpha}\cap K)$. Thus there exists $u_{\alpha}\in U_{\alpha}\cap K$ such that $\Psi(m_{\mu}u_{\alpha}m_{\mu}^{-1})\neq 1$ and 
\begin{equation}
f(m_{\mu})=f(m_{\mu}u_{\alpha})=\overline{\Psi}(m_{\mu}u_{\alpha}m_{\mu}^{-1})f(m_{\mu}).
\end{equation}
Thus $f(m_{\mu})=0$. \EndProof\\

\noindent\textbf{Remark:} The proof of Lemma \ref{CompIndSup} uses the fact that $\Psi$ has conductor $\frak{p}$.\\

Now we can construct a geometric basis for $\mathrm{ind}_{\overline{U}}^{G}(\overline{\Psi})^{K}$. Let $\lambda\in X_{*}(\textbf{A})^{++}$. Define

\begin{equation}
\phi_{\lambda}(g)\stackrel{\mathrm{def}}{=}
\begin{cases}
\delta_{\overline{P}}(m_{\lambda})^{1/2}\overline{\Psi}(u)&,\text{ if } g=um_{\lambda}k\in\overline{U}m_{\lambda}K;\\
0&,\text{ otherwise.}
\end{cases}
\end{equation}
The function $\phi_{\lambda}$ is well-defined because $\lambda$ is strictly dominant.
\begin{lemma}
The set $\{\phi_{\lambda}|\lambda\in X_{*}(\textbf{A})^{++}\}$ is a basis for $\mathrm{ind}_{\overline{U}}^{G}(\overline{\Psi})^{K}$.
\end{lemma}

\textbf{Proof:} This follows from Lemma \ref{CompIndSup} and the Iwasawa decomposition.\EndProof\\

Now we can investigate the $\mathcal{H}_{K}$-module structure of $\mathrm{ind}_{\overline{U}}^{G}(\overline{\Psi})^{K}$ by comparing it to $\mathcal{H}_{IK}$ using the twisted Satake transform. The twisted Satake transform  $\mathcal{S}_{\Psi}^{\prime}:C_{c}(G/I)\rightarrow \mathrm{ind}_{\overline{U}}^{G}(\overline{\Psi})^{I}$ is defined by 
\begin{equation}
\mathcal{S}_{\Psi}^{\prime}(f)(m)=\int_{\overline{U}}f(um)\Psi(u)du,
\end{equation}
where the Haar measure of $\overline{U}$ is normalized so that $\mathrm{meas}(\overline{U}\cap I)=1$. (This is \textit{not} the normalization used to define the Satake transform.) A direct calculation shows that $\mathcal{S}_{\Psi}^{\prime}$ is a homomorphism of left $\mathcal{H}$-modules. Specifically, $\mathcal{S}_{\Psi}^{\prime}(\varrho(f_{1})\cdot f_{2})=\varrho(f_{1})\cdot \mathcal{S}_{\Psi}^{\prime}(f_{2})$, where $f_{1}\in \mathcal{H}$ and $f_{2}\in C_{c}(G/I)$.

We are primarily interested in $\mathcal{S}_{\Psi}\stackrel{\mathrm{def}}{=}\mathcal{S}_{\Psi}^{\prime}|_{\mathcal{H}_{IK}}$. First, we study the image of $\mathcal{S}_{\Psi}$.

\begin{lemma}\label{TwistedOnto}
Let $\lambda \in X_{*}(\textbf{A})^{++}$. Then $\mathcal{S}_{\Psi}(\theta_{\lambda}^{K})=\phi_{\lambda}$. In particular, the twisted Satake transform $\mathcal{S}_{\Psi}:\mathcal{H}_{IK}\rightarrow (\mathrm{ind}_{\overline{U}}^{G}\overline{\Psi})^{K}$ is surjective.
\end{lemma}
\textbf{Proof:} By Lemma \ref{CompIndSup}, it suffices to compute $\mathcal{S}_{\Psi}(\theta_{\lambda}^{K})(m_{\mu})$, where $\mu\in X_{*}(\textbf{A})^{++}$. Since $\lambda$ is dominant we have $\theta_{\lambda}^{K}=\delta_{\overline{P}}^{1/2}(m_{\lambda})1_{Im_{\lambda}K}$. Thus

\begin{equation}\label{TwistedCompEqn1}
\mathcal{S}_{\Psi}(\theta_{\lambda}^{K})(m_{\mu})=\delta_{\overline{P}}^{1/2}(m_{\lambda})\int_{\overline{U}}1_{Im_{\lambda}K}(um_{\mu})\Psi(u)du.
\end{equation}

For $u\in\overline{U}$, we show that $um_{\mu}\in Im_{\lambda}K$ implies that $\mu=\lambda$ and $u\in I\cap \overline{U}$. 

Since $\lambda$ is dominant the Iwahori factorization implies $Im_{\lambda}K=(I\cap \overline{U})m_{\lambda}K$. So, if $um_{\mu}\in Im_{\lambda}K\subseteq (I\cap\overline{U})m_{\lambda}K$, then there exists $u^{\prime}\in (I\cap\overline{U}$) and $k\in K$ such that $m_{\mu}=u^{-1}u^{\prime}m_{\lambda}k$. This implies that $(m_{\mu}^{-1}u^{-1}u^{\prime}m_{\mu})(m_{\mu}^{-1}m_{\lambda})\in \overline{P}\cap K=(\overline{U}\cap K)(M\cap K)$. Thus $m_{\mu}^{-1}m_{\lambda}\in M\cap K$, which implies that $\lambda=\mu$; and $m_{\lambda}^{-1}u^{-1}u^{\prime}m_{\lambda}\in (\overline{U}\cap K)$, which implies that $u^{-1}u^{\prime}\in m_{\lambda}(\overline{U}\cap K)m_{\lambda}^{-1}\subsetneq (\overline{U}\cap I)$. The containment is strict because $\lambda$ is \textit{strictly} dominant. Therefore $u\in I\cap \overline{U}$. 

Thus
\begin{equation}\label{TwistedCompEqn2}
\mathcal{S}_{\Psi}(\theta_{\lambda}^{K})(m_{\lambda})=\delta_{\overline{P}}^{1/2}(m_{\lambda})\int_{I\cap \overline{U}}\Psi(u)du=\delta_{\overline{P}}^{1/2}(m_{\lambda}).
\end{equation}
\EndProof\\
Second, we study $\ker \mathcal{S}_{\Psi}$.
\begin{lemma}\label{TwistedKer}
Let $\alpha\in \Delta$, $s$ the simple reflection corresponding to $\alpha$, and $\iota_{\alpha}=1_{I}+\mathcal{T}_{s}$ (the characteristic function of the parahoric subgroup $I_{\alpha}$). Then
\begin{enumerate}
\item $\mathcal{S}_{\Psi}^{\prime}(\iota_{\alpha})=0$;\label{ParahoricKer}
\item $\mathcal{S}_{\Psi}(\theta_{\lambda}^{K}+\theta_{s(\lambda)}^{K})=0$, for all $\lambda\in X_{*}(\textbf{A})$.\label{AltKer}
\end{enumerate} 
\end{lemma}

\textbf{Proof:} (\ref{ParahoricKer}) We will prove the following claims.

\begin{enumerate}[label=\roman*)]
\item $\mathcal{S}_{\Psi}^{\prime}(\iota_{\alpha})(1)=0$.\label{paraker1}
\item If $\overline{U}_{\lambda,w}\stackrel{\mathrm{def}}{=}\{u\in \overline{U}| um_{\lambda}n_{w}\in I_{\alpha}\}\neq \emptyset$, then $\lambda=0$ and $w=1$ or $s$.\label{paraker2}
\item $\mathcal{S}_{\Psi}^{\prime}(\iota_{\alpha})(m_{\lambda}n_{w})=\int_{\overline{U}_{\lambda,w}}\Psi(u)du$.\label{paraker3}
\end{enumerate}

First we show how these claims imply (\ref{ParahoricKer}). Note that $\mathcal{S}_{\Psi}^{\prime}(\iota_{\alpha})$ is determined by its values on representatives of the double cosets $\overline{U}\backslash G/I\cong\{m_{\lambda},m_{\lambda}n_{w}|\lambda\in X_{*}(\textbf{A}),\,w\in W\}$. Thus it suffices to compute $\mathcal{S}_{\Psi}^{\prime}(\iota_{\alpha})(m_{\lambda}n_{w})$, where $\lambda\in X_{*}(\textbf{A})$ and $w\in W$. By \ref{paraker3}, $\mathcal{S}_{\Psi}^{\prime}(\iota_{\alpha})(m_{\lambda}n_{w})=\int_{\overline{U}_{\lambda,w}}\Psi(u)du$. So if $\mathcal{S}_{\Psi}^{\prime}(\iota_{\alpha})(m_{\lambda}n_{w})\neq 0$, then $\overline{U}_{\lambda,w}\neq \emptyset$. In this case, \ref{paraker2} implies that $\lambda=0$ and $w=1$ or $s$, so $m_{\lambda}n_{w}\in I_{\alpha}$. Therefore $\mathcal{S}_{\Psi}^{\prime}(\iota_{\alpha})(m_{\lambda}n_{w})=\mathcal{S}_{\Psi}^{\prime}(\iota_{\alpha})(1)$, which is zero by \ref{paraker1}. 

Now we prove the claims. By definition $\mathcal{S}_{\Psi}(\iota_{\alpha})(1)=\int_{\overline{U}\cap I_{\alpha}}\Psi(u)du=0$, because $\Psi$ is nontrivial on $\overline{U}\cap I_{\alpha}$. This proves item \ref{paraker1}. Item \ref{paraker3} follows from the definition of $\overline{U}_{\lambda,w}$.

Finally we will prove item \ref{paraker2}. Since $I_{\alpha}\subseteq K$ and $n_{w}\in K$ we see that if $u\in \overline{U}_{\lambda,w}$, then $um_{\lambda}\in \overline{P}\cap K=(\overline{U}\cap K)(M\cap K)$. Thus $u\in \overline{U}\cap K$ and $m_{\lambda}\in M\cap K$, which implies that $m_{\lambda}=1$. Thus it suffices to determine $\overline{U}_{0,w}$.

By Tits \cite{T79} Section 3.4 there is a group scheme $\mathscr{G}$ defined over $\mathcal{O}$ such that $\mathscr{G}\times_{\mathcal{O}}F= \textbf{G}$ and $\mathscr{G}(\mathcal{O})=K$. Let $\underline{\textbf{G}}$ be the group scheme over $\kappa=\mathcal{O}/\frak{p}$ defined by reduction mod $\frak{p}$. Since $u,w\in K=\mathscr{G}(\mathcal{O})$ we apply the canonical homomorphism $\mathcal{O}\rightarrow \kappa$ and work in the group $\underline{\textbf{G}}(\kappa)$. We will suppress the underline to avoid clutter.

Let $\textbf{P}_{\alpha}(\kappa)$ be the parabolic subgroup of $\textbf{G}(\kappa)$ that is the image of $I_{\alpha}$ under the mod $\mathfrak{p}$ reduction map. We have $uw\in \overline{\textbf{U}}(\kappa)w\cap \textbf{P}_{\alpha}(\kappa)$. But $\overline{\textbf{U}}(\kappa)w\subset \overline{\textbf{P}}_{\alpha}(\kappa)w\textbf{U}(\kappa)$, and $\textbf{P}_{\alpha}(\kappa)\subseteq \overline{\textbf{P}}_{\alpha}(\kappa)\textbf{U}(\kappa)$. ($\textbf{P}(\kappa) = \textbf{M}(\kappa)\textbf{U}(\kappa)\subseteq \overline{\textbf{P}}(\kappa)\textbf{U}(\kappa)$, and $\textbf{U}_{\alpha}(\kappa)w_{\alpha}\textbf{P}(\kappa)\subseteq \overline{\textbf{P}}_{\alpha}(\kappa)\textbf{U}(\kappa)$.) The intersection $[\overline{\textbf{P}}_{\alpha}(\kappa)w\textbf{U}(\kappa)]\cap[\overline{\textbf{P}}_{\alpha}(\kappa)\textbf{U}(\kappa)]$ is nonempty if and only if $w=1$ or $s$, by the Bruhat decomposition. This completes the proof of item \ref{paraker2}.

(\ref{AltKer}) By definition 

\begin{equation}\label{AltKerEqn1}
\theta_{\lambda}^{K}+\theta_{s(\lambda)}^{K}=(\theta_{\lambda}+\theta_{s(\lambda)})*1_{K}=\frac{1}{[I_{\alpha}:I]}(\theta_{\lambda}+\theta_{s(\lambda)})*\iota_{\alpha}*1_{K}.
\end{equation}
Since $\iota_{\alpha}=1_{I}+\mathcal{T}_{s}$, Corollary \ref{CommId} implies
\begin{equation}\label{AltKerEqn2}
\mathrm{RHS}(\ref{AltKerEqn1})=\frac{1}{[I_{\alpha}:I]}\iota_{\alpha}*(\theta_{\lambda}+\theta_{s(\lambda)})*1_{K}=\frac{1}{[I_{\alpha}:I]}\iota_{\alpha}*(\theta_{\lambda}^{K}+\theta_{s(\lambda)}^{K}).
\end{equation}
We apply the $\mathcal{H}$-module homomorhpism $\mathcal{S}_{\Psi}$ to equations (\ref{AltKerEqn1}) and (\ref{AltKerEqn2}) to get
\begin{equation}
\mathcal{S}_{\Psi}(\theta_{\lambda}^{K}+\theta_{s(\lambda)}^{K})=\frac{1}{[I_{\alpha}:I]}\mathcal{S}_{\Psi}^{\prime}(\iota_{\alpha})*(\theta_{\lambda}^{K}+\theta_{s(\lambda)}^{K}),
\end{equation}
which is zero by item (\ref{ParahoricKer}).
\EndProof\\

\textbf{Remarks:} \textbf{a)} We provide another proof of claim \ref{paraker2} above, when $\lambda=0$. Specifically,  if $u\in \overline{U}\cap K$, $w\in W$, and $uw\in I_{\alpha}$, then $w=1$ or $s$. This approach will utilize the building $\mathscr{B}$. Let $\mathscr{F}_{\alpha}$ be the facet associated to $I_{\alpha}$. Then we know that $I_{\alpha}$ fixes $\mathscr{F}_{\alpha}$ pointwise. Since $w^{-1}u^{-1}\in I_{\alpha}$ it must fix $\mathscr{F}_{\alpha}$. We also know that $w$ maps $\mathscr{A}$ into itself. Thus $u^{-1}\mathscr{F}_{\alpha}$ must be in $\mathscr{A}$. Since $u^{-1}$ acts via an isometry on $\mathscr{B}$ and $u^{-1}$ fixes $-\mathscr{C}$ we claim that $u^{-1}$ must fix $\mathscr{F}_{\alpha}$. To see this we will use the following fact. Let $\mathcal{E}$ be a Euclidean space with distance function $d_{\mathcal{E}}$. For any $y\in \mathcal{E}$ and $d\in \mathbb{R}_{\geq0}$ let $S(y,d)=\{x\in \mathcal{E}|d_{\mathcal{E}}(x,y)=d\}$ be the sphere centered at $y$ or radius $d$. Given $x_{0}\in \mathcal{E}$ and $Y\subset \mathcal{E}$ a set with nonempty interior, we have $\cap_{y\in Y}S(y,d_{\mathcal{E}}(x_{0},y))=\{x_{0}\}$. The apartment $\mathscr{A}$ is a Euclidean space containing the set $-\mathcal{C}$, which has a nonempty interior. If $x_{0}\in \mathscr{F}_{\alpha}$ then $u^{-1}x_{0}\in \cap_{y\in -\mathcal{C}}S(y,d_{\mathcal{E}}(x_{0},y))=\{x_{0}\}$, since for all $y\in -\mathscr{C}$, we have $u^{-1}y=y$. This implies that $w$ fixes $\mathscr{F}_{\alpha}$. Thus $w=1,s$.\\
 
\noindent\textbf{b)} Our proof of Lemma \ref{TwistedKer}, specifically the use of the parahoric subgroup $I_{\alpha}$, forces our choice of a character $\Psi$ of conductor $\mathfrak{p}$. It seems that a different integral operator in place of $\mathcal{S}_{\Psi}$ could be used to circumvent this restriction. Unfortunately this introduces other complications and will not be pursued here. \\

Lemmas \ref{TwistedOnto} and \ref{TwistedKer} directly imply the following corollary.

\begin{corollary}\label{TwistedKerCor}
The kernel of the $\mathcal{H}_{K}$-module homomorphism $\mathcal{S}_{\Psi}:\mathcal{H}_{IK}\rightarrow \mathrm{ind}_{\overline{U}}^{G}(\overline{\Psi})^{K}$ is 
\begin{equation}
\ker(\mathcal{S}_{\Psi})= \mathrm{span}\{\theta_{\lambda}^{K}-(-1)^{\ell(w)}\theta_{w\lambda}^{K}|\lambda\in X_{*}(\textbf{A}),\, w\in W\}.
\end{equation}
\end{corollary}

Now we introduce the fundamental diagram of $\mathcal{H}_{K}$-modules.
\begin{equation}\label{MainDiagram}
\begin{tikzcd}
\mathcal{H}_{IK}\arrow[d,"\mathcal{S}"]\arrow[r, "{\mathcal{S}_{\Psi}}"]& \mathrm{ind}_{\overline{U}}^{G}(\overline{\Psi})^{K}\arrow[d,dashed,"{j}"]\\
\mathbb{C}[X_{*}(\textbf{A})]\arrow[r,"\text{alt}"]&\mathrm{alt}(\mathbb{C}[X_{*}(\textbf{A})])
\end{tikzcd}
\end{equation}
In the next proposition we define $j$ and show that it is an isomorphism of right $\mathcal{H}_{K}$-modules.

\begin{proposition}\label{indIso}
$\mathcal{S}(\ker(\mathcal{S}_{\Psi}))=\ker(\mathrm{alt})$. Thus we can define the $\mathcal{H}_{K}$-module isomorphism $j:\mathrm{ind}_{\overline{U}}^{G}(\overline{\Psi})^{K}\rightarrow \mathrm{alt}(\mathbb{C}[X_{*}(\textbf{A})])$ by the formula $j(\phi)=\mathrm{alt}(S(f))$, where $f\in\mathcal{H}_{IK}$ such that $\mathcal{S}_{\Psi}(f)=\phi$. In particular, for $\mu\in X_{*}(\textbf{A})^{++}$ we have $j(\phi_{\mu})=\mathrm{alt}(e^{\mu})$.
\end{proposition}

\textbf{Proof:} This follows from Theorem \ref{SavinIso}, Lemma \ref{TwistedOnto}, and Corollary \ref{TwistedKerCor}.\EndProof

\section{Unramified Whittaker Function}\label{UWF}

Recall that $R=\mathbb{C}[X_{*}(\textbf{A})]$ and $S\supseteq R$ is a $\mathbb{C}$-algebra that is an integral domain. In this section we prove that the $S$-valued Whittaker function $\mathcal{W}$ (introduced in Subsection \ref{WF}) satisfies a family of recursions, Proposition \ref{WhittakerRecur}.

Let $\langle-,-\rangle=\langle-,-\rangle_{S,\Psi}:\mathrm{ind}_{\overline{U}}^{G}(\overline{\Psi})\otimes_{\mathbb{C}}\mathrm{Ind}_{\overline{U}}^{G}(S_{\Psi}) \rightarrow S$ be the bilinear form defined by $\langle\phi_{1},\phi_{2}\rangle=\int_{\overline{U}\backslash G}\phi_{1}(g)\phi_{2}(g)dg$, where the measure on $\overline{U}\backslash G$ is induced from the Haar measure on $G$ such that $\mathrm{meas}(K)=1$. Note that this pairing is $G$-invariant.

Now we establish a few basic properties of $\langle-,\mathcal{W}\rangle$. 

\begin{lemma}\label{WhitEval}
Let $\lambda\in X_{*}(\textbf{A})^{++}$ and let $f\in \mathrm{Ind}_{\overline{U}}^{G}(S_{\Psi})^{K}$. Then
\begin{equation}
\langle\phi_{\lambda},f\rangle=\delta_{\overline{P}}^{-1/2}(m_{\lambda})f(m_{\lambda}).
\end{equation}
\end{lemma}

\textbf{Proof:} Let $F\in C_{c}^{\infty}(G/K,S)$ such that $\int_{\overline{U}}F(ug)du=\phi_{\lambda}(g)f(g)$. Then
\begin{equation}\label{WhitEvalEqn1}
\langle\phi_{\lambda},f\rangle\stackrel{\mathrm{def}}{=}\int_{\overline{U}\backslash G}\phi_{\lambda}(g)f(g)dg =\int_{G}F(g)dg.
\end{equation}
By the Iwasawa decomposition and the right $K$-invariance of $F$,
\begin{equation}\label{WhitEvalEqn2}
\mathrm{RHS}(\ref{WhitEvalEqn1})=\int_{\overline{P}}\int_{K}F(pk)dkd_{\ell}p=\int_{\overline{P}}F(p)d_{\ell}p.
\end{equation}
By equation (9) in Cartier \cite{C79}, page 145, 
\begin{equation}\label{WhitEvalEqn3}
\mathrm{RHS}(\ref{WhitEvalEqn2})=\int_{M}\int_{\overline{U}}\delta_{\overline{P}}(m)^{-1} F(um)dudm.
\end{equation}
By integrating over $\overline{U}$ and applying the definition of $\phi_{\lambda}$ we have
\begin{equation}\label{WhitEvalEqn4}
\mathrm{RHS}(\ref{WhitEvalEqn3})=\int_{M}\delta_{\overline{P}}(m)^{-1}\phi_{\lambda}(m)f(m)dm=\delta_{\overline{P}}^{-1/2}(m_{\lambda})f(m_{\lambda}).
\end{equation}
\EndProof\\
 
\begin{lemma}\label{WhitPairing}
Let $\phi\in \mathrm{ind}_{\overline{U}}^{G}(\overline{\Psi})^{K}$ and $f\in \mathcal{H}_{K}$. Then 
\begin{equation}
\langle\varrho(f)\cdot \phi,\mathcal{W}\rangle=\mathcal{S}(\hat{f})\langle\phi,\mathcal{W}\rangle.
\end{equation}
\end{lemma}

\textbf{Proof:} The $G$-invariance of the pairing implies that $\langle\varrho(f)\cdot \phi,\mathcal{W}\rangle=\langle\phi,\varrho(\hat{f})\cdot\mathcal{W}\rangle$. The result follows from the identity $\varrho(\hat{f})\cdot \mathcal{W}=\mathcal{S}(\hat{f})\mathcal{W}$.\EndProof

\begin{lemma} \label{HKModHom}
Let $\lambda\in X_{*}(\textbf{A})^{+}$ and let $f\in \mathcal{H}_{IK}$. Then
\begin{equation}\label{HKMod1}
\langle\mathcal{S}_{\Psi}(f*A_{\lambda}),\mathcal{W}\rangle=\mathrm{ch}V_{\lambda}\langle\mathcal{S}_{\Psi}(f),\mathcal{W}\rangle.
\end{equation}

In particular, if $\mu\in X_{*}(\textbf{A})^{++}$ and $f=\theta_{\mu}^{K}$, then 
\begin{equation}\label{HKMod2}
\langle\mathcal{S}_{\Psi}(\theta_{\mu}^{K}*A_{\lambda}),\mathcal{W}\rangle=\delta_{\overline{P}}^{-1/2}(m_{\mu})\mathrm{ch}V_{\lambda}\mathcal{W}(m_{\mu})
\end{equation}
\end{lemma}

\textbf{Proof:} The identities $\mathcal{S}_{\Psi}(f*A_{\lambda})=\mathcal{S}_{\Psi}(f)*A_{\lambda}$ and $f*A_{\lambda}=\varrho(\hat{A}_{\lambda})\cdot f$, and Lemma \ref{WhitPairing} imply Equation (\ref{HKMod1}). Equation (\ref{HKMod1}) and Lemma \ref{WhitEval} imply Equation (\ref{HKMod2}).\EndProof\\

Let $\lambda\in X_{*}(\textbf{A})^{+}$ and $\mu\in X_{*}(\textbf{A})^{++}$. Now we will describe an explicit formula for $\phi_{\mu}*A_{\lambda}$ in terms of characters of finite dimensional representations of ${^{L}\textbf{G}}^{\dagger}$.

\begin{lemma} \label{TensorDecomp}
Let $\lambda\in X_{*}(\textbf{A})^{+}$ and $\mu\in X_{*}(\textbf{A})^{++}$. Let $c_{\mu,\lambda}^{\eta}\in \mathbb{C}$ be defined as in subsection \ref{DualGroup}. Then

\begin{equation}\label{TensorDecompEqn}
\phi_{\mu}*A_{\lambda}=\sum_{\eta\in X_{*}(\textbf{A})^{++}}c_{\mu,\lambda}^{\eta}\phi_{\eta}.
\end{equation}
\end{lemma}

\textbf{Proof:} Since $j$ is an isomorphism and $\mathrm{alt}\circ \mathcal{S}=j\circ \mathcal{S}_{\Psi}$ it suffices to show that 
\begin{equation*}
\mathrm{alt}\circ \mathcal{S}(\theta_{\mu}^{K}*A_{\lambda})=\mathrm{alt}\circ \mathcal{S}(\sum_{\eta}c_{\mu,\lambda}^{\eta}\theta_{\eta}^{K}).
\end{equation*} 

By Theorem \ref{SavinIso}, we have
\begin{equation}\label{TensorDecompEqn1}
\mathrm{alt}\circ \mathcal{S}(\theta_{\mu}^{K}*A_{\lambda})=\mathrm{ch}V_{\lambda}\cdot\mathrm{alt}(e^{\mu}).
\end{equation}

We multiply and divide by $\mathrm{alt}(e^{\rho^{\vee}})$ and apply the Weyl character formula with respect to the group ${^{L}\textbf{G}}^{\dagger}$, specifically $\mathrm{alt}(e^{\mu})=\mathrm{ch}V_{\mu-\rho^{\vee}} \mathrm{alt}(e^{\rho^{\vee}})$, (temporarily working in the field of fractions of $\mathbb{C}[\mathscr{X}]$) to get 

\begin{equation}\label{TensorDecompEqn2}
\mathrm{RHS}(\ref{TensorDecompEqn1})=\mathrm{ch}V_{\mu-\rho^{\vee}}\cdot\mathrm{ch}V_{\lambda}\cdot\mathrm{alt}(e^{\rho^{\vee}})
\end{equation}
Next we apply the identity $\mathrm{ch}V_{\mu-\rho^{\vee}}\cdot \mathrm{ch}V_{\lambda}=\sum_{\eta}c_{\mu,\lambda}^{\eta}\mathrm{ch}V_{\eta-\rho^{\vee}}$ followed by the Weyl character formula to get
\begin{align}
\mathrm{RHS}(\ref{TensorDecompEqn2})=&\sum_{\eta}c_{\mu,\lambda}^{\eta}\mathrm{ch}V_{\eta-\rho^{\vee}}\cdot\mathrm{alt}(e^{\rho^{\vee}})\nonumber\\
=&\sum_{\eta}c_{\mu,\lambda}^{\eta}\mathrm{alt}(e^{\eta})\label{TensorDecompEqn3}.
\end{align}

Again by Theorem \ref{SavinIso}, we have
\begin{equation}\label{TensorDecompEqn4}
\mathrm{RHS}(\ref{TensorDecompEqn3})=\mathrm{alt}\circ \mathcal{S}(\sum_{\eta}c_{\mu,\lambda}^{\eta}\theta_{\eta}^{K}).
\end{equation}
\EndProof\\

\noindent\textbf{Remark:} If $\rho^{\vee}\in X_{*}(\textbf{A})$, then Lemma \ref{TensorDecomp} implies $\mathrm{ind}_{\overline{U}}^{G}(\Psi)^{K}\cong \mathcal{H}_{K}\cdot \phi_{\rho^{\vee}}$.\\

We conclude this section with a family of recursions for $\mathcal{W}$.

\begin{proposition}\label{WhittakerRecur}
Let $\lambda\in X_{*}(\textbf{A})^{+}$ and $\mu\in X_{*}(\textbf{A})^{++}$. Define $c_{\mu,\lambda}^{\eta}\in\mathbb{C}$ be such that $\mathrm{ch}V_{\mu-\rho^{\vee}}\cdot\mathrm{ch}V_{\lambda}=\sum_{\eta}c_{\mu,\lambda}^{\eta}\mathrm{ch}V_{\eta-\rho^{\vee}}$. (Recall Subsection \ref{DualGroup}.) Then 
\begin{equation}
\delta_{\overline{P}}^{-1/2}(m_{\mu})\mathrm{ch}V_{\lambda}\cdot\mathcal{W}(m_{\mu})=\sum_{\eta\in X_{*}(\textbf{A})^{++}}c_{\mu,\lambda}^{\eta}\delta_{\overline{P}}^{-1/2}(m_{\eta})\mathcal{W}(m_{\eta}).
\end{equation}
\end{proposition}

\textbf{Proof:} Apply $\langle-,\mathcal{W}\rangle$ to equation (\ref{TensorDecompEqn}) and use Lemmas \ref{WhitEval} and \ref{HKModHom}.\EndProof\\

In the next two sections we prove two special cases of the Casselman-Shalika formula based on Proposition \ref{WhittakerRecur}. In Section \ref{CSAdjointSec} we prove the Casselman-Shalika formula where $\textbf{G}$ is a semisimple group of adjoint type and $\Psi$ has conductor $\mathfrak{p}$ or $\mathcal{O}$ (defined in Section \ref{CSAdjointSec}). In Section \ref{CSFormulaSec} we prove the Casselman-Shalika formula where $\textbf{G}$ is an arbitrary connected unramified group, but $\Psi$ is of conductor $\mathfrak{p}$.


\section{Casselman-Shalika Formula for Adjoint Groups}\label{CSAdjointSec}

In this section we suppose that $\textbf{G}$ is semisimple of adjoint type and prove the Casselman-Shalika formula for $R$-valued spherical Whittaker functions associated to characters of conductor $\mathfrak{p}$ or $\mathcal{O}$. Note that $\textbf{G}$ adjoint implies that $\rho^{\vee}\in X_{*}(\textbf{A})$.

\begin{theorem}[Conductor $\mathfrak{p}$]\label{CSAdjoint}
Let $\lambda\in X_{*}(\textbf{A})^{+}$. Then
\begin{equation}
\mathcal{W}(m_{\lambda+\rho^{\vee}})=\delta_{\overline{P}}^{1/2}(m_{\lambda})\mathrm{ch}V_{\lambda}\cdot\mathcal{W}(m_{\rho^{\vee}}).
\end{equation}
\end{theorem}

\textbf{Proof:} Take $\mu=\rho^{\vee}$ in Proposition \ref{WhittakerRecur}. Then 
\begin{equation}\label{AdRecur}
\mathrm{ch}V_{\mu-\rho^{\vee}}\cdot \mathrm{ch}V_{\lambda}=\sum_{\eta}c_{\mu,\lambda}^{\eta}\mathrm{ch}V_{\eta-\rho^{\vee}}=\mathrm{ch}V_{\lambda}.
\end{equation}
Thus
\begin{equation}
\delta_{\overline{P}}^{-1/2}(m_{\rho^{\vee}})\mathrm{ch}V_{\lambda}\cdot\mathcal{W}(m_{\rho^{\vee}})=\delta_{\overline{P}}^{-1/2}(m_{\lambda+\rho^{\vee}})\mathcal{W}(m_{\lambda+\rho^{\vee}}),
\end{equation}
from which it follows that
\begin{equation}
\mathcal{W}(m_{\lambda+\rho^{\vee}})=\delta_{\overline{P}}^{1/2}(m_{\lambda})\mathrm{ch}V_{\lambda}\cdot\mathcal{W}(m_{\rho^{\vee}}).
\end{equation}
\EndProof\\

\noindent\textbf{Remarks:} a) The proof in this section is valid for any $\textbf{G}$ such that $\rho^{\vee}\in X_{*}(\textbf{A})$.\\

\noindent b) We take this opportunity to correct the formula appearing in Theorem 6.1 in \cite{N13}. It is off by a factor of $\delta_{P}^{1/2}(m_{\rho^{\vee}})$.\\

Now we treat the case where $\Psi$ has conductor $\mathcal{O}$. We accomplish this by relating the Whittaker function of conductor $\mathcal{O}$ to a Whittaker function of conductor $\mathfrak{p}$. 

\begin{lemma}\label{CondSwap}
Suppose that $\Psi_{\mathcal{O}}$ is a non-degenerate character of $\overline{U}$ of conductor $\mathcal{O}$, $\mathfrak{W}_{\mathcal{O}}$ is a Whittaker functional with respect to $\Psi_{\mathcal{O}}$, and $\mathcal{W}_{\mathcal{O}}=\mathfrak{W}_{\mathcal{O}}(gv_{0})$. Then
\begin{enumerate}
\item $\Psi_{\mathfrak{p}}(u)\stackrel{\mathrm{def}}{=}\Psi_{\mathcal{O}}(m_{-\rho^{\vee}}um_{-\rho^{\vee}}^{-1})$ is a non-degenerate character of conductor $\mathfrak{p}$;
\item $\mathfrak{W}_{\mathfrak{p}}\stackrel{\mathrm{def}}{=}\mathfrak{W}_{\mathcal{O}}\circ\varrho(m_{-\rho^{\vee}})$ is a Whittaker functional with respect to $\Psi_{\mathfrak{p}}$;
\item $\mathcal{W}_{\mathfrak{p}}(g)\stackrel{\mathrm{def}}{=}\mathcal{W}_{\mathcal{O}}(m_{-\rho^{\vee}}g)$ is the spherical Whittaker function associated to $\mathfrak{W}_{\mathfrak{p}}$.
\end{enumerate} 
\end{lemma}

\textbf{Proof:} This follows directly from the definitions.\EndProof

\begin{proposition}[Conductor $\mathcal{O}$] \label{CSAdO}
Let $\mathcal{W}_{\mathcal{O}}(g)$ be a spherical Whittaker function  for the semisimple adjoint group $G$ of conductor $\mathcal{O}$. Let $\lambda\in X_{*}(\textbf{A})^{+}$. Then 
\begin{equation}
\mathcal{W}_{\mathcal{O}}(m_{\lambda})=\delta_{\overline{P}}^{1/2}(m_{\lambda})\mathrm{ch}V_{\lambda}\cdot\mathcal{W}_{\mathcal{O}}(1)
\end{equation}
\end{proposition}
\textbf{Proof:} This follows directly from Lemma \ref{CondSwap} and Theorem \ref{CSAdjoint}.\EndProof\\

\textbf{Remark:} This argument can be applied to arbitrary characters $\Psi$.


\section{Casselman-Shalika Formula Via Recursion}\label{CSFormulaSec}

In this section we present a proof of the Casselman-Shalika formula for $R$-valued spherical Whittaker functions of conductor $\mathfrak{p}$ where $\textbf{G}$ is a connected unramified group. The main obstruction that prohibits the argument of Section \ref{CSAdjointSec} is that $\rho^{\vee}$ may not be an element of $X_{*}(\textbf{A})$. 

Let $Q$ be the field of fractions of $R$ and define the $Q$-vector space $V^{\prime}=\oplus_{\mu\in X_{*}(\textbf{A})^{++}}Qe_{\mu}$ with standard basis elements $e_{\mu}$. Recall the coefficients $c_{\lambda,\mu}^{\eta}$ defined in Subsection \ref{DualGroup}. Let 
\begin{equation*}
V^{\prime\prime}\stackrel{\mathrm{def}}{=}\mathrm{span}_{Q}(\mathrm{ch}V_{\lambda}e_{\mu}-\sum c_{\lambda,\mu}^{\eta}e_{\eta}|\lambda\in X_{*}(\textbf{A})^{+},\mu\in X_{*}(\textbf{A})^{++}),
\end{equation*}
and let $p:V^{\prime}\rightarrow V\stackrel{\mathrm{def}}{=}V^{\prime}/V^{\prime\prime}$ be the canonical quotient map.

We can define two functionals as follows. Let $\mathscr{W}:V^{\prime}\rightarrow Q$ by $\mathscr{W}(e_{\mu})=\delta(m_{\mu})^{-1/2}\mathcal{W}(m_{\mu})$, and let  $\omega:V^{\prime}\rightarrow Q$ defined by $e_{\mu}\mapsto \mathrm{alt}(e^{\rho^{\vee}})\cdot\mathrm{ch}V_{\mu-\rho^{\vee}}=\mathrm{alt}(e^{\mu})$ (Weyl character formula). (Note that $\mathrm{ch}V_{\mu-\rho^{\vee}}$ is a character of ${^{L}\textbf{G}}^{\dagger}$.)

\begin{lemma}
The functionals $\omega,\mathscr{W}$ factor through $V$.  Moreover, $V\neq 0$.
\end{lemma} 

\textbf{Proof:} Proposition \ref{WhittakerRecur} implies that $\mathscr{W}$ factors through $V$. The Weyl character formula implies that $\omega$ is well-defined and factors through $V$. Note that $V\neq 0$ because $\omega\neq 0$.\EndProof\\

 If we can prove that $V\cong Q$ as a $Q$-vector space, then $\omega$ and $\mathscr{W}$ are proportional and the Casselman-Shalika formula follows.

\begin{proposition}\label{dimV}
The $Q$-linear map $\omega$ defines an isomorphism $V\cong Q$.
\end{proposition}

\textbf{Proof:} Let $\mu\in X_{*}(\textbf{A})^{++}$. Then $\omega(e_{\mu})=\mathrm{alt}(e^{\mu})\neq 0$. So it suffices to show that $\ker\omega=V^{\prime\prime}$. 

The multiplication map $R\otimes_{\mathbb{C}}R\rightarrow R$ restricts to give the following the exact sequence of $R^{W}$-modules
\begin{equation}\label{multExact}
0\rightarrow \ker(\mathrm{mult})\rightarrow R^{W}\otimes_{\mathbb{C}}\mathrm{alt}(R)\stackrel{\mathrm{mult}}{\rightarrow} \mathrm{alt}(R)\rightarrow 0.
\end{equation}
Since $\mathrm{alt}(R)=\mathbb{C}[X_{*}(\textbf{A})+\rho^{\vee}]^{W}\mathrm{alt}(e^{\rho^{\vee}})$, the definition of the ring structure of $\mathrm{Rep}({^{L}\textbf{G}}^{\dagger})$ implies that 
\begin{multline}
\ker(\mathrm{mult})=\\\mathrm{span}_{\mathbb{C}}(\mathrm{ch}V_{\lambda}\otimes \mathrm{ch}V_{\mu-\rho^{\vee}}\mathrm{alt}(e^{\rho^{\vee}})-\sum c_{\lambda,\mu}^{\eta}(1\otimes \mathrm{ch}V_{\eta-\rho^{\vee}}\mathrm{alt}(e^{\rho^{\vee}}))|\lambda\in X_{*}(\textbf{A})^{+}, \mu\in X_{*}(\textbf{A})^{++}).
\end{multline}
Consider the functor $Q\otimes_{R^{W}}-$. It is exact because it is the composition of the exact functors $Q^{W}\otimes_{R^{W}}-$ and $Q\otimes_{Q^{W}}-$. The first functor $Q^{W}\otimes_{R^{W}}-$ is exact because $Q^{W}$ is the field of fractions of $R^{W}$; the second functor $Q\otimes_{Q^{W}}-$ is exact because $Q$ is a free $Q^{W}$-module since $Q^{W}$ is a field. Thus if we apply $Q\otimes_{R^{W}}-$ to the exact sequence (\ref{multExact}) we get the exact sequence
\begin{equation}
0\rightarrow Q\otimes_{R^{W}}\ker(\mathrm{mult})\rightarrow Q\otimes_{\mathbb{C}}\mathrm{alt}(R)\stackrel{\mathrm{mult}}{\rightarrow} Q\rightarrow 0.
\end{equation}
Let $\Xi:Q\otimes_{\mathbb{C}}\mathrm{alt}(R)\rightarrow V^{\prime}$ be the $Q$-module isomorphism defined by $1\otimes \mathrm{alt}(e^{\mu})\mapsto e_{\mu}$. By definition we have the following commutative diagram.

\begin{equation}\label{MultComm}
\begin{tikzcd}
Q\otimes_{\mathbb{C}}\mathrm{alt}(R)\arrow[d, hook, two heads,"\Xi"]\arrow[r,"{\mathrm{mult}}"]& Q\arrow[d,equals]\\
V^{\prime}\arrow[r,"\omega"]&Q
\end{tikzcd}
\end{equation}

Thus $\ker\omega= \Xi(Q\otimes_{R^{W}}\ker(\mathrm{mult}))=V^{\prime\prime}$.\EndProof

\begin{theorem}\label{CSFormulaFinal} Let $\mu\in X_{*}(\textbf{A})^{++}$, then 
\begin{equation}\label{CSFormulaEqn}
\mathcal{W}(m_{\mu})=r\delta_{\overline{P}}^{1/2}(m_{\mu})\mathrm{ch}V_{\mu-\rho^{\vee}}\cdot\mathrm{alt}(e^{\rho^{\vee}}),
\end{equation} 
where $r\in R$ is a normalization factor depending on the choice of Whittaker functional $\mathfrak{W}$.
\end{theorem}

\textbf{Proof:} By Proposition \ref{dimV} we have that $\mathrm{dim}_{Q}(V)=1$, thus $\mathrm{Hom}_{Q}(V,Q)\cong Q\omega$. Since $\mathscr{W}\in \mathrm{Hom}_{Q}(V,Q)$ there exists $r\in Q$ such that $\mathscr{W}=r\omega$. By the definition of $\mathscr{W}$, equation (\ref{CSFormulaEqn}) follows. In fact, $r\in R$ because $\mathcal{W}$ is $R$-valued, $R$ is a unique factorization domain, and $r$ is independent of $\mu$.\EndProof


\section{The General Case}\label{CSGeneral}

In this section we deduce the Casselman-Shalika formula for a general connected unramified group $\textbf{G}^{\prime}$ and $\Psi$ a character of conductor $\mathfrak{p}$. We use the notation of Subsection \ref{RUPS}. In particular, $\textbf{G}=\textbf{G}^{\prime}/Z(\textbf{G}^{\prime})$, $\textbf{T}=\textbf{G}^{\prime}/\mathcal{D}\textbf{G}^{\prime}$, $\textbf{G}^{\prime\prime}=\textbf{G}\times_{\mathrm{Spec}(F)}\textbf{T}$, $R^{\prime}=\mathbb{C}[X_{*}(\textbf{A}^{\prime})]$, and $R^{\prime\prime}=\mathbb{C}[X_{*}(\textbf{A}^{\prime\prime})]$.

Let $\mathfrak{W}^{\prime}\in \mathrm{Hom}_{(R^{\prime},\overline{U})}(i_{G^{\prime},\overline{P}^{\prime}}((\chi_{\mathrm{univ}}^{\prime})^{-1}),R^{\prime}_{\Psi})$ be an $R^{\prime}$-valued Whittaker functional for the group $G^{\prime}$. This gives rise to the $R^{\prime}$-valued spherical Whittaker function $\mathcal{W}^{\prime}(g^{\prime})=\mathfrak{W}^{\prime}(g^{\prime}v_{0}^{\prime})$. We will relate this Whittaker function to a Whittaker function on the semisimple group $G$ of adjoint type. Since $R^{\prime}\hookrightarrow R^{\prime\prime}$ (recall line (\ref{cocharinj})), $\mathfrak{W}^{\prime}$ and $\mathcal{W}^{\prime}$ may also be viewed as $R^{\prime\prime}$-valued.

Lemma \ref{GroupTransfer2} and Lemma \ref{WhittFunctExtScale} imply that 
\begin{equation}
\mathrm{Hom}_{(R^{\prime},\overline{U})}(i_{G^{\prime},\overline{P}^{\prime}}((\chi_{\mathrm{univ}}^{\prime})^{-1}),R^{\prime\prime}_{\Psi})\cong\mathrm{Hom}_{(R^{\prime\prime},\overline{U}^{\prime})}(i_{G^{\prime\prime},\overline{P}^{\prime\prime}}((\chi_{\mathrm{univ}}^{\prime\prime})^{-1}),R^{\prime\prime}_{\Psi})
\end{equation}
as $R^{\prime\prime}$-modules. Thus $\mathfrak{W}^{\prime}$ corresponds to $\mathfrak{W}^{\prime\prime}\in \mathrm{Hom}_{(R^{\prime\prime},\overline{U}^{\prime})}(i_{G^{\prime\prime},\overline{P}^{\prime\prime}}((\chi_{\mathrm{univ}}^{\prime\prime})^{-1}),R^{\prime\prime}_{\Psi})$ such that $\mathfrak{W}^{\prime}=\mathfrak{W}^{\prime\prime}|_{1\otimes i_{G^{\prime},\overline{P}^{\prime}}((\chi_{\mathrm{univ}}^{\prime})^{-1})}$, where $1\otimes i_{G^{\prime},\overline{P}^{\prime}}((\chi_{\mathrm{univ}}^{\prime})^{-1})\subseteq i_{G^{\prime\prime},\overline{P}^{\prime\prime}}((\chi_{\mathrm{univ}}^{\prime\prime})^{-1})$ according to Lemma \ref{GroupTransfer2}.

This relationship between Whittaker functionals induces a relationship between spherical Whittaker functions. Let $\mathcal{W}^{\prime\prime}((g,t))=\mathfrak{W}^{\prime\prime}((g,t)v_{0}^{\prime\prime})$. 
\begin{proposition}\label{toProd} Let $g^{\prime}\in G^{\prime}$. Then
\begin{equation}
\mathcal{W}^{\prime}(g^{\prime})=\mathcal{W}^{\prime\prime}(\pi(g^{\prime})).
\end{equation}
\end{proposition}

\textbf{Proof:} Under the isomorphism of Lemma \ref{GroupTransfer2}, $v_{0}^{\prime\prime}=1\otimes v_{0}^{\prime}$. Thus if we write $\pi(g^{\prime})=(g,t)$ then
\begin{equation*}
\mathcal{W}^{\prime}(g^{\prime})=\mathfrak{W}^{\prime\prime}(1\otimes g^{\prime}v_{0}^{\prime})
=\mathfrak{W}^{\prime\prime}(\pi(g^{\prime})v_{0}^{\prime\prime})
=\mathcal{W}^{\prime\prime}(\pi(g^{\prime})).
\end{equation*}
\EndProof\\

With Proposition \ref{toProd} we reduce the computation of a spherical Whittaker function on $G^{\prime}$ to one on $G^{\prime\prime}$; next we reduce the computation of a spherical Whittaker function on $G^{\prime\prime}$ to one on $G$. 

We begin with some notation and a relationship between the Whittaker functionals. Let $\mathfrak{W}^{\prime\prime}\in \mathrm{Hom}_{(R^{\prime\prime},\overline{U}^{\prime})}(i_{G^{\prime\prime},\overline{P}^{\prime\prime}}((\chi_{\mathrm{univ}}^{\prime\prime})^{-1}),R^{\prime\prime}_{\Psi})$. Lemma \ref{ProdPSeries} and Lemma \ref{WhittFunctExtScale} imply that 
\begin{equation}
\mathrm{Hom}_{(R^{\prime\prime},\overline{U})}(i_{G^{\prime\prime},\overline{P}^{\prime\prime}}((\chi_{\mathrm{univ}}^{\prime\prime})^{-1}),R^{\prime\prime}_{\Psi})\cong\mathrm{Hom}_{(R,\overline{U})}(i_{G,\overline{P}}(\chi_{\mathrm{univ}}^{-1}),R^{\prime\prime}_{\Psi})
\end{equation}
as $R^{\prime\prime}$-modules. Thus $\mathfrak{W}^{\prime\prime}$ corresponds to a $\mathfrak{W}\in \mathrm{Hom}_{(R,\overline{U})}(i_{G,\overline{P}}(\chi_{\mathrm{univ}}^{-1}),R^{\prime\prime}_{\Psi})$ such that $\mathfrak{W}=\mathfrak{W}^{\prime\prime}|_{1\otimes i_{G,\overline{P}}(\chi_{\mathrm{univ}}^{-1})}$.

Again the relationship between Whittaker functionals induces a relationship between spherical Whittaker functions. Let $\mathcal{W}(g)=\mathfrak{W}(gv_{0})$. Recall that $\xi$ is the tautological character of $T$ (Subsection \ref{RUPS}).

\begin{proposition}\label{toAd} Let $(g,t)\in G^{\prime\prime}$. Then
\begin{equation}
\mathcal{W}^{\prime\prime}((g,t))=\xi(t)^{-1}\mathcal{W}(g).
\end{equation}

\end{proposition}

\textbf{Proof:} Under the isomorphism of Lemma \ref{ProdPSeries}, $v_{0}^{\prime\prime}=1\otimes v_{0}$. Thus
\begin{equation*}
\mathcal{W}^{\prime\prime}((g,t))=\xi(t)^{-1}\mathfrak{W}^{\prime\prime}(1\otimes gv_{0})
=\xi(t)^{-1}\mathfrak{W}(gv_{0})
=\xi(t)^{-1}\mathcal{W}(g).
\end{equation*} \EndProof\\

We combine Proposition \ref{toProd} and Proposition \ref{toAd} to express $\mathcal{W}^{\prime}$ in terms of $\mathcal{W}$, a spherical Whittaker function of a semisimple adjoint group $G$.

\begin{proposition}\label{WhittRed} Let $g^{\prime}\in G^{\prime}$, $g\in G$, and $t\in T$ such that $\pi(g^{\prime})=(g,t)$. Then
\begin{equation}
\mathcal{W}^{\prime}(g^{\prime})=\xi(t)^{-1}\mathcal{W}(g).
\end{equation}

\end{proposition}

The next theorem is the Casselman-Shalika formula for a general connected reductive group $G^{\prime}$.

\begin{theorem}\label{CSP}
Let $\mu^{\prime}\in X_{*}(A^{\prime})^{++}$ and let $\pi_{*}(\mu^{\prime})=(\mu,\lambda)\in X_{*}(\textbf{A})\times X_{*}(\textbf{T})=X_{*}(\textbf{A}^{\prime\prime})$. Then
\begin{equation}\label{CSPEqn}
\mathcal{W}^{\prime}(m_{\mu^{\prime}})=\xi^{-1}(t_{\lambda})\delta_{\overline{P}}^{1/2}(m_{\mu})\mathrm{ch}V_{\mu-\rho^{\vee}}\cdot\mathcal{W}(m_{\rho^{\vee}}).
\end{equation}
\end{theorem}

\textbf{Proof:} This is a direct consequence of Proposition \ref{WhittRed} and Theorem \ref{CSAdjoint}.\EndProof\\

\textbf{Remark:} In Theorem \ref{CSP}, $\mathrm{ch}V_{\mu-\rho^{\vee}}$ is a character of the Langlands dual group of $\textbf{G}^{\prime}/Z(\textbf{G}^{\prime})$ and thus lives in $R^{\prime\prime}$. We also have $\mathcal{W}(\rho^{\vee})\in R^{\prime\prime}$ and $\xi^{-1}(t_{\lambda})\in R^{\prime\prime}$. However, since $\mathcal{W}^{\prime}$ is $R^{\prime}$-valued $\mathrm{RHS}(\ref{CSPEqn})$ is in $R^{\prime}$.


\section{Conductor $\mathcal{O}$}\label{ConductorO}

Our focus on Whittaker functions of conductor $\mathfrak{p}$ is a by-product of our proof technique. However, for global applications one must study Whittaker functions of conductor $\mathcal{O}$. This means that for all $\alpha\in \Delta$, $\Psi_{\alpha}$ is trivial on $U_{-\alpha}\cap K$ and nontrivial on any larger subgroup of $U_{-\alpha}$. We accomplish this by relating a Whittaker function of conductor $\mathcal{O}$ to a Whittaker function of conductor $\mathfrak{p}$, for which Theorem \ref{CSP} provides a formula. 

We will retain the notation from Section \ref{CSGeneral} and augment the notation of $\Psi$, $\mathfrak{W}$, and $\mathcal{W}$ to indicate the conductor. For example, $\Psi_{\mathcal{O}}$ will be a character of $\overline{U}$ of conductor $\mathcal{O}$ and $\Psi_{\mathfrak{p}}$ be a character of conductor $\mathfrak{p}$.

All of the pieces are already in place because the results of Section \ref{CSGeneral} with the exception of Theorem \ref{CSP} are independent of the conductor of $\Psi$. (The conductor $\mathfrak{p}$ property is used in our proof of Lemma \ref{TwistedKer}, an important ingredient for Theorem \ref{CSP}.) 

\begin{theorem}\label{CSO}
Let $\lambda^{\prime}\in X_{*}(A^{\prime})$ be dominant. Let $\lambda\in X_{*}(\textbf{A})$ be dominant, and $\mu\in X_{*}(T)$  such that $\pi_{*}(\lambda^{\prime})=(\lambda,\mu)$. Then
\begin{equation}
\mathcal{W}^{\prime}_{\mathcal{O}}(m_{\lambda^{\prime}})=\xi^{-1}(t_{\mu})\delta_{\overline{P}}^{1/2}(m_{\lambda})\mathrm{ch}V_{\lambda}\cdot\mathcal{W}_{\mathcal{O}}^{\prime}(1).
\end{equation}
\end{theorem}

\textbf{Proof:} Apply Propositions \ref{WhittRed}, \ref{CSAdO}, and Theorem  \ref{CSP} and note that $\mathcal{W}_{\mathcal{O}}^{\prime}(1)=\mathcal{W}_{\mathcal{O}}(1)$. \EndProof \\

\noindent\textbf{Remarks:} 
\begin{enumerate}
\item $\xi^{-1}(t_{\mu})$ and $\mathrm{ch}V_{\lambda}$ need not be in the image of $R^{\prime}\hookrightarrow R^{\prime\prime}=R\otimes_{\mathbb{C}}\mathbb{C}[T/T^{\circ}]$ individually, but the product $\xi^{-1}(t_{\mu})\mathrm{ch}V_{\lambda}\cdot\mathcal{W}_{\mathcal{O}}^{\prime}(1)$ is.
\item One can apply the argument of this section to spherical Whittaker functions associated to characters $\Psi$ with arbitrary conductor.
\end{enumerate}


\end{document}